%% file: main.tex
\documentclass[journal]{new-aiaa}
\usepackage[utf8]{inputenc}
\usepackage{textcomp}

\usepackage{graphicx}
\usepackage{subcaption}
\DeclareCaptionLabelFormat{aiaa}{#2)}
\usepackage{amsmath}
\usepackage[version=4]{mhchem}
\usepackage{siunitx}
\usepackage{longtable,tabularx}
\usepackage[printonlyused]{acronym}
\input{sections/__acronyms}
\hypersetup{hidelinks}
\usepackage{dsfont}
\usepackage{todonotes}
\usepackage[table]{xcolor}

\title{Theoretical and Computational Advances in Birkhoff Pseudospectral Transcription for Aerospace Optimal Control\footnote{Presented at the AIAA SCITECH 2027 Forum, Orlando, FL, January 11--15, 2027.}}

\author{Riccardo Minnozzi \footnote{PhD Student, Department of Industrial Engineering} and Marco Sagliano \footnote{Full Professor, Department of Industrial Engineering, AIAA Associate Fellow}}
\affil{University of Bologna, Via L. Montaspro 97, Forlì (FC), Italy, 47121}

\begin{document}

\maketitle

\input{sections/abstract}

\input{sections/introduction}
\input{sections/birkhoff_theory}
\input{sections/differential_constraints}
\input{sections/numerical_advancements}
\input{sections/results}
\input{sections/conclusions}

\section*{Acknowledgments}
\label{sec:acknowledgments}
Funded by the European Union (\ac{erc}, STARGATE, 101169959). Views and opinions expressed are however those of the author(s) only and do not necessarily reflect those of the European Union or the European Research Council. Neither the European Union nor the granting authority can be held responsible for them.
\begin{figure}[h!]
\centering
\includegraphics[width=3.25in]{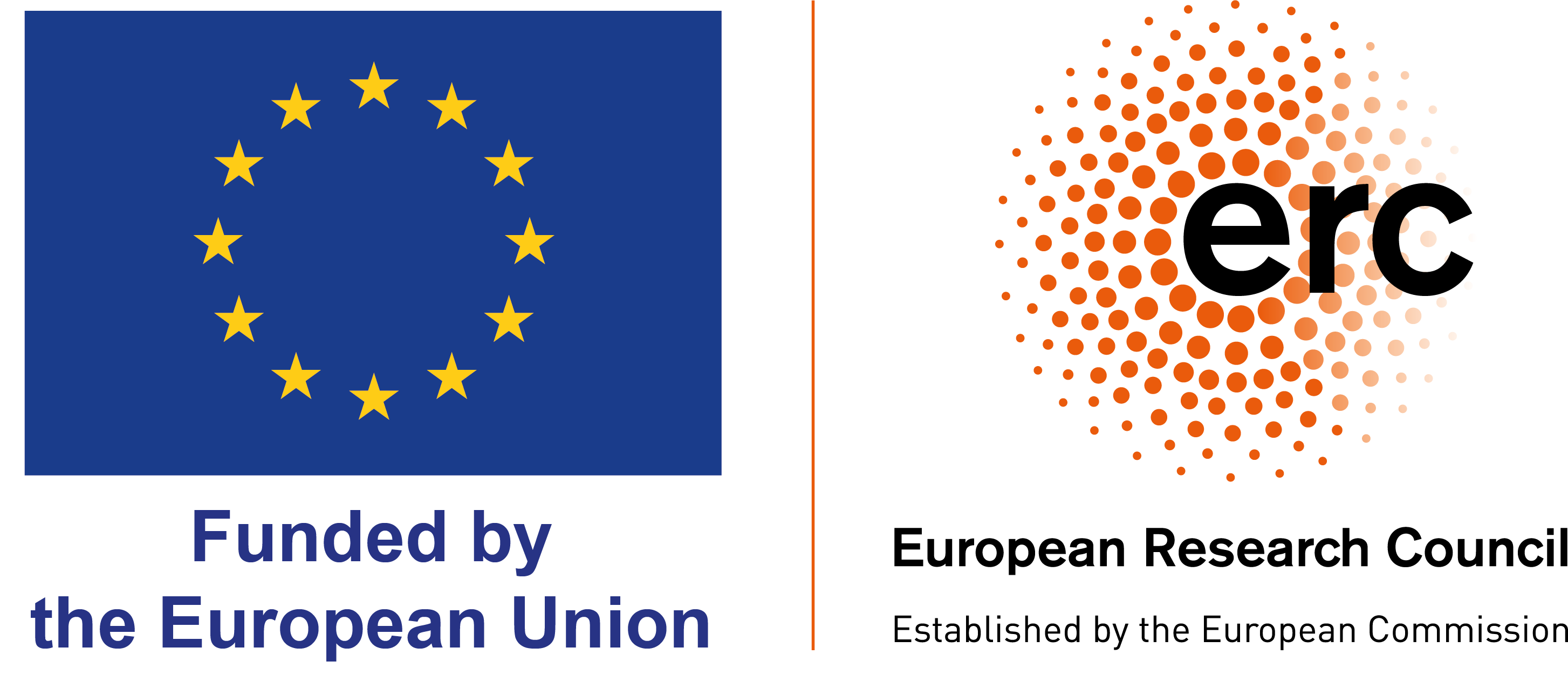}
\captionsetup{labelformat=empty}
\caption{}
\end{figure}

\bibliography{references}

\end{document}

%% file: sections/__acronyms.tex
\acrodef{nlp}[NLP]{Nonlinear Program}
\acrodef{raptor}[RAPTOR]{Rapid Pseudospectral Trajectory Optimization and Robustification}
\acrodef{ocp}[OCP]{Optimal Control Problem}
\acrodef{lgr}[LGR]{Legendre-Gauss-Radau}
\acrodef{lgl}[LGL]{Legendre-Gauss-Lobatto}
\acrodef{flgr}[fLGR]{flipped-Legendre-Gauss-Radau}
\acrodef{rlv}[RLV]{Reusable Launch Vehicle}
\acrodef{ps}[PS]{Pseudospectral}
\acrodef{ipopt}[IPOPT]{Interior Point OPTimizer}
\acrodef{snopt}[SNOPT]{Sparse Nonlinear OPTimizer}
\acrodef{worhp}[WORHP]{We Optimize Really Huge Problems}
\acrodef{rtls}[RTLS]{Return To Launch Site}
\acrodef{kkt}[KKT]{Karush-Kuhn-Tucker}
\acrodef{erc}[ERC]{European Research Council}
\acrodef{stargate}[STARGATE]{STochastic Algorithm for Robust Guidance Analysis and Trajectory Estimation}
\acrodef{jrn}[JRN]{Jacobian Row Normalization}
\acrodef{pjrn}[PJRN]{Projected Jacobian Row Normalization}

%% file: sections/abstract.tex
\begin{abstract}
\par Birkhoff pseudospectral methods replace the differentiation matrices of classical transcriptions with integration operators whose conditioning does not deteriorate with polynomial order. This paper presents theoretical and computational advancements toward their practical use in aerospace optimal control. The initial-point-anchored Birkhoff integration operator is constructed explicitly for flipped Legendre-Gauss-Radau grids and related exactly to the corresponding Lagrange differentiation matrix and, through a reflection, to the standard Radau Birkhoff operator, whose conditioning properties it inherits. Rate constraints on nonlinear functions of states and controls are enforced directly through the pseudospectral operators and, for control-rate limits, shown to be primal and dual equivalent to state augmentation. Non-affine time mappings concentrate the nodes around localized trajectory features, and composing affine scaling and equilibration transformations provides a systematic way to construct scaling strategies. The methods are implemented in the open-source MATLAB framework RAPTOR and tested across a range of challenging optimal control problems, including a dense-grid orbit-raising problem, a rate-limited landing scenario, a Space Shuttle entry problem, and a five-phase reusable-launch-vehicle mission.
\end{abstract}

\acresetall

%% file: sections/introduction.tex
\section{Introduction}
\label{sec:introduction}

\par Trajectory optimization is central to the design and guidance of aerospace vehicles across a wide range of missions, including interplanetary transfers, atmospheric entry, powered descent and landing, and the ascent and recovery of reusable launch vehicles. Such trajectory-design problems are naturally formulated as \acp{ocp} \cite{bryson2018_AppliedOptimalControl}, whose numerical approaches are commonly classified as indirect or direct methods. Indirect methods apply the necessary conditions of optimality to reformulate the \ac{ocp} as a multi-point boundary-value problem \cite{longuski2014_OptimalControlAerospace}, whereas direct methods discretize the continuous problem and transcribe it into a finite-dimensional parametric optimization problem \cite{betts2010_PracticalMethodsOptimal}.

\par Direct methods have become widely adopted for aerospace trajectory optimization thanks to their robustness in enforcing nonlinear dynamics, operational constraints, and complex mission architectures \cite{betts2010_PracticalMethodsOptimal, ross2012_ReviewPseudospectralOptimal, malyuta2022_ConvexOptimizationTrajectory}. Among them, \ac{ps} methods approximate the continuous state trajectories using high-order polynomial representations and enforce the system dynamics at a set of nonuniform collocation nodes, thereby providing algebraic approximations of the differential and integral operators \cite{fornberg2005_PracticalGuidePseudospectral}. For sufficiently smooth solutions, these methods exhibit spectral convergence with increasing polynomial order, meaning that complex trajectories can be accurately captured with low-order discretizations \cite{ross2012_ReviewPseudospectralOptimal}. Moreover, appropriate covector-mapping theorems enable the reconstruction of the continuous-time costates from the multipliers of the discretized problem \cite{ross2001_PseudospectralTransformationCovectors, gong2008_ConnectionsCovectorMapping}. Different \ac{ps} formulations arise from both the selected node distribution and the interpolation representation of the states, with classical methods commonly based on Lagrange interpolation at Gauss-type nodes \cite{fornberg2005_PracticalGuidePseudospectral}.

\par Aerospace applications of \ac{ps} methods range from low-thrust trajectory optimization \cite{ross2015_PrimerPontryaginsPrinciple, ross2007_LowThrustHighAccuracyTrajectory, kelly2021_GeostationaryDebrisMitigation} to atmospheric entry, powered descent, and reusable-launch-vehicle trajectory design \cite{huneker2015_SPARTANImprovedGlobal, sagliano2017_RadauPseudospectralMethod, sagliano2024_SixDegreeofFreedomRocketLanding, sagliano2024_SixDegreesofFreedomAeroPropulsiveEntry, sagliano2025_SixDegreeofFreedomEntryPowered}. Beyond their traditional coupling with \ac{nlp} solvers, \ac{ps} discretizations have also been incorporated into convex trajectory-optimization frameworks \cite{sagliano2019_GeneralizedHpPseudospectralConvex}. Despite their favorable convergence properties, high-order Lagrange \ac{ps} discretizations exhibit an important numerical limitation: the conditioning of differentiation-matrix-based formulations deteriorates approximately as $\mathcal{O}(N^2)$ with the number of discretization nodes $N$ \cite{koeppen2019_FastMeshRefinement}. Multi-interval and $hp$-adaptive approaches mitigate this limitation by distributing the approximation among multiple lower-order polynomial segments; nevertheless, the deteriorating conditioning of the underlying differentiation formulation limits the practical use of very high-order discretizations.

\par Birkhoff-based \ac{ps} formulations provide an alternative approach to this problem. By constructing the state approximation from mixed boundary-value and derivative information, the differential transcription can be recast in an integration-matrix form whose conditioning can remain essentially independent of the discretization order for relevant boundary-condition configurations \cite{proulx2024_ImplementationsUniversalBirkhoff, ross2024_UniversalBirkhoffTheory}. An extensive theoretical framework on universal Birkhoff methods is introduced in Refs.~\cite{ross2023_UniversalBirkhoffPseudospectral, ross2024_UniversalBirkhoffTheory}, covering the construction and covector mapping of Birkhoff \ac{ps} methods built on Legendre and Chebyshev grids for Gauss, Radau, and Lobatto node families. These works, however, do not explicitly address \ac{flgr} grids, which are commonly adopted for aerospace problems with prescribed initial conditions \cite{sagliano2017_RadauPseudospectralMethod}.

\par Beyond favorable operator conditioning, the practical solution of aerospace \acp{ocp} requires additional considerations on the formulation of constraints, the resolution of the discretization, and the scaling of the resulting problem. Rate constraints are commonly enforced through state augmentation \cite{wang2017_ConstrainedTrajectoryOptimization, wang2020_ImprovedSequentialConvex}, at the cost of additional states and controls. In differentiation-matrix-based \ac{ps} methods, these can be imposed directly \cite{nie2021_EfficientImplementationRate}, but, to the best of the authors' knowledge, neither an equivalent formulation for integration-based \ac{ps} methods nor its relation to the dual structure of state augmentation has been established. Localized features such as control switches are typically resolved by increasing the polynomial order or subdividing the domain \cite{gong2008_SpectralAlgorithmPseudospectral, sagliano2019_GeneralizedHpPseudospectralConvex}, while non-affine time mappings \cite{ross2007_LowThrustHighAccuracyTrajectory, ross2020_EnhancementsDIDOOptimal} have received limited attention for trajectory resolution purposes. Finally, the conditioning of the \ac{ps} operators does not necessarily extend to the complete \ac{nlp}, therefore effective scaling strategies are required \cite{betts2010_PracticalMethodsOptimal, ross2018_ScalingBalancingHighPerformance}.

\par To address these aspects, this work introduces theoretical and computational advancements toward the practical use of Birkhoff \ac{ps} methods for aerospace \acp{ocp}. The theoretical contributions are the following. First, the Birkhoff-$a$ integration operator is constructed explicitly for \ac{flgr} grids and related to both the corresponding Lagrange differentiation matrix \cite{garg2010_UnifiedFrameworkNumerical} and the \ac{lgr} Birkhoff-$b$ operator of the universal Birkhoff theory \cite{ross2023_UniversalBirkhoffPseudospectral, ross2024_UniversalBirkhoffTheory}, from which it inherits its conditioning properties. Second, differential constraints on nonlinear functions of states and controls are enforced directly through the \ac{ps} operators, and their primal and dual equivalence with state augmentation is shown for control-rate limits. On the computational side, non-affine time mappings \cite{ross2007_LowThrustHighAccuracyTrajectory, ross2020_EnhancementsDIDOOptimal} are exploited to concentrate the \ac{ps} nodes around localized trajectory features, and the composition of affine scaling and equilibration transformations is introduced as a systematic approach to construct \ac{ocp} scaling strategies, recovering \ac{pjrn} \cite{sagliano2017_RadauPseudospectralMethod} as a particular case. These developments are implemented within the \ac{raptor}\footnote{solver: \url{https://github.com/raptor-solver/raptor-matlab.git}, \newline docs: \url{https://raptor-solver.github.io/}} framework, a modular MATLAB transcription engine for multiphase \acp{ocp} developed within the context of the \ac{erc}-funded \ac{stargate} project \cite{minnozzi2026_STARGATERobustGuidance}, and demonstrated on four benchmark problems, the last of which showcases their integration within a branching multiphase \ac{ocp}.

\par The remainder of this paper is organized as follows. Section~\ref{sec:birkhoff-theory} summarizes \ac{ps} integration methods and illustrates the Birkhoff formulation for \ac{flgr} grids, together with its connections to the Lagrange formulation and to standard Radau operators. Section~\ref{sec:differential-constraints} describes direct differential-constraint enforcement and its equivalence to state augmentation. Section~\ref{sec:numerical-advancements} focuses on the numerical advances that led to the unified transcription framework implemented in \ac{raptor}, including non-affine time mappings and the composition of affine scaling techniques. Section~\ref{sec:results} presents the corresponding landing, orbit-raising, atmospheric-entry, and multiphase \ac{rlv} demonstrations.

%% file: sections/birkhoff_theory.tex
\section{Birkhoff Theory for Pseudospectral Optimal Control}
\label{sec:birkhoff-theory}

\par In this section we summarize the underlying theory on differentiation-based \ac{ps} methods, present the Birkhoff interpolants used for the \ac{ps} formulations in Refs.~\cite{ross2019_FurtherResultsFast, ross2020_EnhancementsDIDOOptimal, ross2023_UniversalBirkhoffPseudospectral, ross2024_UniversalBirkhoffTheory}, and derive their formulation when applied to \ac{flgr} grids. We start by introducing the trajectories for states $\mathbf{x}(\tau)\in\mathbb{R}^{n_x}$, controls $\mathbf{u}(\tau)\in\mathbb{R}^{n_u}$, and a vector of static parameters $\mathbf{p}\in\mathbb{R}^{n_p}$. The dynamic variables are defined on the computational domain $\tau\in[-1,1]$, also referred to as pseudotime. Physical time is determined by the mapping $t=\mathcal{T}(\tau;t_0,t_f,\mathbf{p})$, where $t_0$ and $t_f$ are the phase endpoints and need not be fixed. The corresponding time-dilation factor is
\begin{equation}
    \frac{\mathrm{d}t}{\mathrm{d}\tau}
    =
    \rho(\tau;t_0,t_f,\mathbf{p})
    =
    \frac{\partial\mathcal{T}}{\partial\tau}(\tau;t_0,t_f,\mathbf{p})
    >0
    \label{eq:time-dilation-definition}
\end{equation}
The shortened notation $\rho(\tau)$ is adopted for convenience when dependence on $t_0$, $t_f$, and $\mathbf{p}$ is understood. A more detailed discussion on time mappings is provided in Sec.~\ref{subsec:non-affine-time-mappings}. The time-dilation factor maps the state dynamics $\mathbf{f}(\mathbf{x}, \mathbf{u}, t, \mathbf{p})$ from physical time to pseudotime as
\begin{equation}
    \frac{\mathrm{d}\mathbf{x}}{\mathrm{d}\tau}=\rho(\tau)\mathbf{f}\!\left(\mathbf{x}(\tau),\mathbf{u}(\tau),t(\tau),\mathbf{p}\right)
    \label{eq:pseudotime-dynamics}
\end{equation}

\par \ac{ps} methods use a global interpolant to approximate the full state trajectory: given a set of $N_d$ discretization nodes $\{\tau_i^d\}_{i=1}^{N_d}$, a Lagrange interpolant $\mathbf{x}^L(\tau)$ provides the lowest-order polynomial fit satisfying all nodal values, and it is defined as
\begin{equation}
    \mathbf{x}^{L}(\tau)=\sum_{i=1}^{N_d}\mathbf{X}_i\,\ell_i(\tau),\qquad \text{with } \mathbf{X}_i\approx\mathbf{x}(\tau_i^d)
    \label{eq:lagrange-state-interpolation}
\end{equation}
where $\ell_i(\tau)$ is the cardinal Lagrange polynomial
\begin{equation}
    \ell_i(\tau)=\prod_{\substack{k=1\\k\neq i}}^{N_d}\frac{\tau-\tau_k^d}{\tau_i^d-\tau_k^d},\qquad\text{with } \ell_i(\tau_k^d)=\delta_{ik}
    \label{eq:lagrange-cardinal-polynomials}
\end{equation}
and $\delta_{ik}$ is the Kronecker delta.

\par The traditional approach to \ac{ps} methods for optimal control is the definition of a discrete differentiation operator, hence why in this work we equivalently use the terms "Lagrange formulation" and "differentiation-based \ac{ps} methods" \cite{fornberg2005_PracticalGuidePseudospectral, sagliano2017_RadauPseudospectralMethod, ross2012_ReviewPseudospectralOptimal}. The differentiation operator maps the discretized quantities onto a collocation grid $\{\tau_j^c\}_{j=1}^{N_c}$ (which need not be identical to the discretization grid) by differentiating Eq.~\eqref{eq:lagrange-state-interpolation} directly, resulting in the generally rectangular differentiation matrix $\mathbf{D}\in\mathbb{R}^{N_c\times N_d}$
\begin{equation}
    \left.\frac{\mathrm{d}\mathbf{x}^{L}}{\mathrm{d}\tau}\right|_{\tau_j^c}=\sum_{i=1}^{N_d}D_{ji}\mathbf{X}_i\qquad D_{ji}=\ell_i'(\tau_j^c)
    \label{eq:lagrange-differentiation}
\end{equation}
Collecting the dynamic variables via their collocated values $\mathbf{X}_j^c=\mathbf{x}^{L}(\tau_j^c)$, $\mathbf{U}_j\approx\mathbf{u}(\tau_j^c)$, and $\rho_j=\rho(\tau_j^c)$, the Lagrange dynamics defect constraints are imposed through 
\begin{equation}
    \sum_{i=1}^{N_d}D_{ji}\mathbf{X}_i-\rho_j\mathbf{f}\left(\mathbf{X}_j^c,\mathbf{U}_j,t_j,\mathbf{p}\right)=\mathbf{0}\qquad j=1,\ldots,N_c
    \label{eq:lagrange-defects}
\end{equation}

\par The defect constraint in Eq.~\eqref{eq:lagrange-defects} explicitly contains the differentiation matrix $\mathbf{D}$. For the grids commonly used in \ac{ps} methods for optimal control, the condition number of this operator increases as $\mathcal{O}(N^2)$ with increasing polynomial order \cite{koeppen2019_FastMeshRefinement}. 
Since the operator enters the Jacobian of the corresponding \ac{nlp}, the numerical performance of the associated optimization algorithms may deteriorate significantly \cite{betts2010_PracticalMethodsOptimal}. Multi-interval $hp$-\ac{ps} and knotting methods provide one strategy to avoid very high-order global differentiation matrices \cite{gong2008_SpectralAlgorithmPseudospectral, sagliano2019_GeneralizedHpPseudospectralConvex}. Birkhoff methods, on the other hand, modify the interpolation representation itself so that high-order discretizations can be constructed while mitigating the conditioning growth associated with the \ac{ps}-specific portions of the transcribed problem.

\subsection{Overview of Birkhoff Methods}
\label{subsec:birkhoff-overview}

\par We now summarize the relevant theory on Birkhoff methods for optimal control, as it constitutes the main framework for \ac{ps} integration methods. While Lagrange interpolants approximate a function from a set of values at the nodes, Birkhoff interpolants can be built from a combination of function and derivative values \cite{wang2014_WellConditionedCollocationMethod, ross2019_FurtherResultsFast}. For the definition of a \ac{ps} integration method, one boundary value and a collocation of derivative values are sufficient: while multiple different formulations can be constructed depending on the choice of boundary value, a natural choice is the initial condition $\mathbf{X}_a$, and we refer to the corresponding formulation as the Birkhoff-$a$ version \cite{ross2023_UniversalBirkhoffPseudospectral, ross2024_UniversalBirkhoffTheory}. Introducing the virtual derivative variables
\begin{equation}
    \mathbf{V}_j\approx\left.\frac{\mathrm{d}\mathbf{x}}{\mathrm{d}\tau}\right|_{\tau_j^c}\qquad j=1,\ldots,N_c
    \label{eq:birkhoff-derivative-variables}
\end{equation}
the Birkhoff-$a$ interpolant for the state trajectory can be defined as $\mathbf{x}^{B}(\tau)$
\begin{equation}
    \mathbf{x}^{B}(\tau)=\mathbf{X}_a+\sum_{j=1}^{N_c}\mathbf{V}_j B_j^a(\tau)\qquad \mathbf{X}_a\approx\mathbf{x}(-1)
    \label{eq:birkhoff-interpolant}
\end{equation}
and it exhibits the cardinality property
\begin{equation}
    B_j^a(-1)=0\qquad \left(B_j^a\right)'(\tau_k^c)=\delta_{jk}
    \label{eq:birkhoff-cardinality}
\end{equation}

\par The Birkhoff theory developed in Refs.~\cite{ross2019_FurtherResultsFast, ross2023_UniversalBirkhoffPseudospectral, ross2024_UniversalBirkhoffTheory} constructs the Birkhoff interpolants by assuming that the state derivative can be represented via a Lagrange polynomial (as it was done for the state trajectory in Eq.~\eqref{eq:lagrange-state-interpolation}). The basis functions for the Birkhoff interpolant are therefore obtained via integration of the cardinal Lagrange polynomials $\ell_j^c$ associated with the collocation grid, defined as in Eq.~\eqref{eq:lagrange-cardinal-polynomials} with the $N_c$ collocation nodes $\tau_j^c$ in place of the $N_d$ discretization nodes $\tau_i^d$
\begin{equation}
    B_j^a(\tau)=\int_{-1}^{\tau}\ell_j^c(s)\,\mathrm{d}s
    \label{eq:birkhoff-basis-integral}
\end{equation}

\par The Birkhoff integration matrix is obtained by evaluating Eq.~\eqref{eq:birkhoff-basis-integral} at the collocation nodes $B_{ij}=B_j^a(\tau_i^c)$. The Birkhoff dynamic defects are imposed through the reconstruction equation
\begin{equation}
    \mathbf{X}_i^c=\mathbf{X}_a+\sum_{j=1}^{N_c}B_{ij}\mathbf{V}_j\qquad i=1,\ldots,N_c
    \label{eq:birkhoff-reconstruction}
\end{equation}
Depending on the chosen collocation grid, the final state may not be included (as it is the case for \ac{lgr} grids), and therefore does not enter Eq.~\eqref{eq:birkhoff-reconstruction}. To maintain a consistent definition of the full state history (i.e., including potentially non-collocated boundary values $\mathbf{X}_b$), Refs.~\cite{ross2023_UniversalBirkhoffPseudospectral, ross2024_UniversalBirkhoffTheory} introduce a calculus-theorem constraint that binds the endpoint state values with the state derivative via quadrature
\begin{equation}
    \mathbf{X}_b-\mathbf{X}_a=\sum_{j=1}^{N_c}w_j^{B}\mathbf{V}_j\qquad w_j^{B}=B_j^a(+1)=\int_{-1}^{1}\ell_j^c(s)\,\mathrm{d}s
    \label{eq:birkhoff-calculus}
\end{equation}
To complete the Birkhoff dynamics defects, the virtual derivative variables $\mathbf{V}_j$ are connected to the physical-time dynamics through
\begin{equation}
    \mathbf{V}_j=\rho_j\mathbf{f}\left(\mathbf{X}_j^c,\mathbf{U}_j,t_j,\mathbf{p}\right)\qquad j=1,\ldots,N_c
    \label{eq:birkhoff-dynamics}
\end{equation}

\par Equations~\eqref{eq:birkhoff-reconstruction} and \eqref{eq:birkhoff-dynamics} replace the differentiation defects of Eq.~\eqref{eq:lagrange-defects} with a linear integration relation and a set of local nonlinear dynamics equations. The resulting Jacobian contains blocks of the form $[\mathbf{B} |  -\mathbf{I}_{N_c}]$, instead of the high-order differentiation matrix $\mathbf{D}$. The conditioning of this discretization-specific structure grows as $\mathcal{O}(\sqrt{N_c})$ \cite{koeppen2019_FastMeshRefinement}, and it can behave as well as $\mathcal{O}(1)$ if prescribed boundary conditions are exploited \cite{proulx2024_ImplementationsUniversalBirkhoff, ross2024_UniversalBirkhoffTheory}.

\subsection{Birkhoff Interpolants on Flipped-Radau Grids}
\label{subsec:new-birkhoff-interpolants}

\par The construction of the Birkhoff integration operator depends on the grid choice for the specific problem: while previous works provide formulations for Gauss, Radau, and Lobatto grids based on Chebyshev and Legendre polynomials \cite{ross2023_UniversalBirkhoffPseudospectral, ross2024_UniversalBirkhoffTheory}, this section complements existing work by deriving an explicit operator for the \ac{flgr} grid on the initial-endpoint-anchored Birkhoff-$a$ formulation, and connects it to the regular \ac{lgr} operators on the Birkhoff-$b$ version in Sec.~\ref{subsubsec:flgr-birkhoff-reflection}.

\par The \ac{flgr} grid contains the final endpoint $\tau=1$, and it is defined by the roots of the polynomial $P_M(\tau)-P_{M-1}(\tau)$, where $P_M$ is the $M$-th order Legendre polynomial \cite{fornberg2005_PracticalGuidePseudospectral} and $M=N_c$ is the number of collocation points. The Birkhoff-$a$ version is particularly natural in this setting, as the initial state value $\mathbf{X}_a \approx \mathbf{x}(\tau=-1)$ is not included in the collocation grid. Since our implementation adopts interchangeable transcriptions through both Lagrange and Birkhoff interpolants, we adopt the same notation as the Lagrange formulation in Eq.~\eqref{eq:lagrange-state-interpolation} and define the discretization grid as $\boldsymbol{\tau}^{d}=[-1,\tau_1^c,\ldots,\tau_M^c]^{T}$. Thus, the state is represented at $M+1$ discretization nodes, while its derivative and the dynamics are collocated at the $M$ \ac{flgr} nodes. As in standard Radau \ac{ps} methods, the controls are defined at the collocation nodes, and the transcription does not determine their value at the noncollocated node $\tau=-1$. In \ac{raptor}, an additional control value $\mathbf{U}_a\approx\mathbf{u}(-1)$ is associated with this node: unless it enters a constraint (see Sec.~\ref{sec:differential-constraints}), its value is obtained by extrapolation of the collocated controls.

\par To construct the Birkhoff integration matrix directly, we consider the cardinal Lagrange polynomials $\ell_j^c$ associated with the \ac{flgr} collocation grid, introduced in Eq.~\eqref{eq:birkhoff-basis-integral}. The polynomial has degree at most $M-1$ and can therefore be expanded via a linear combination of the Legendre polynomial basis $P_0,\ldots,P_{M-1}$. The normalization of this basis is \cite{fornberg2005_PracticalGuidePseudospectral}
\begin{equation}
    \gamma_k
    =
    \int_{-1}^{1}P_k^2(s)\,\mathrm{d}s
    =
    \frac{2}{2k+1}
    \qquad k=0,\ldots,M-1
    \label{eq:legendre-normalization}
\end{equation}
The \ac{flgr} quadrature weights $w_j^{\mathrm{fLGR}}$ compute the exact quadrature integration of all polynomials $q(s)$ of degree at most $2M-2$, meaning
\begin{equation}
    \int_{-1}^{1}q(s)\,\mathrm{d}s
    =
    \sum_{r=1}^{M}w_r^{\mathrm{fLGR}}q(\tau_r^c)
    \label{eq:flgr-quadrature-rule}
\end{equation}

\par The coefficient of the $k$-th order term $P_k$ in the Legendre expansion of $\ell_j^c$ is
\begin{equation}
    \frac{1}{\gamma_k}
    \int_{-1}^{1}\ell_j^c(s)P_k(s)\,\mathrm{d}s
    =
    \frac{w_j^{\mathrm{fLGR}}P_k(\tau_j^c)}{\gamma_k}
    \label{eq:flgr-cardinal-legendre-coefficient}
\end{equation}
This equality follows from Eq.~\eqref{eq:flgr-quadrature-rule}: the integrand has degree at most $2M-2$, and the cardinality property in Eq.~\eqref{eq:lagrange-cardinal-polynomials} reduces the quadrature sum to its $j$-th term. Hence,
\begin{equation}
    \ell_j^c(\tau)
    =
    w_j^{\mathrm{fLGR}}
    \sum_{k=0}^{M-1}
    \frac{P_k(\tau_j^c)P_k(\tau)}{\gamma_k}
    \label{eq:flgr-cardinal-legendre-expansion}
\end{equation}

\par To evaluate Eq.~\eqref{eq:birkhoff-basis-integral}, we define the integrated Legendre modes $\mathcal{I}_k(\tau)$ as
\begin{equation}
    \mathcal{I}_k(\tau)
    =
    \int_{-1}^{\tau}P_k(s)\,\mathrm{d}s
    =
    \begin{cases}
        1+\tau & k=0\\
        \dfrac{P_{k+1}(\tau)-P_{k-1}(\tau)}{2k+1} & k\geq1
    \end{cases}
    \label{eq:integrated-legendre}
\end{equation}
Since $P_k(\tau)$ is the only $\tau$-dependent factor in Eq.~\eqref{eq:flgr-cardinal-legendre-expansion}, integrating it term by term from $-1$ to $\tau_i^c$ replaces each $P_k$ with $\mathcal{I}_k$ and directly gives the entries of the pseudospectral integration matrix $\mathbf{B}\in\mathbb{R}^{M\times M}$:
\begin{equation}
    B_{ij}
    =
    w_j^{\mathrm{fLGR}}
    \sum_{k=0}^{M-1}
    \frac{P_k(\tau_j^c)\mathcal{I}_k(\tau_i^c)}{\gamma_k}
    \qquad
    i,j=1,\ldots,M
    \label{eq:flgr-birkhoff-matrix}
\end{equation}

\par Because the final \ac{flgr} node is $\tau_M^c=+1$, the last row of $\mathbf{B}$ integrates the cardinal polynomials over the full interval $\tau\in[-1, 1]$, thus 
\begin{equation}
    B_{Mj} = \int_{-1}^1\ell_j^c(s)\mathrm{d}s=w_j^{\mathrm{fLGR}}
    \label{eq:birkhoff-flgr-final-row}
\end{equation}
For this reason, the final row of Eq.~\eqref{eq:birkhoff-reconstruction} for the specific case of \ac{flgr} grids (where $\mathbf{X}_M^c=\mathbf{X}_b$) is the discrete calculus theorem of Eq.~\eqref{eq:birkhoff-calculus}, meaning that the separate enforcement of Eq.~\eqref{eq:birkhoff-calculus} would be a duplicate constraint in a Birkhoff-$a$ \ac{flgr} formulation, and it is therefore removed.

\par The state reconstruction constraint from Eq.~\eqref{eq:birkhoff-reconstruction} is expressed in matrix form as
\begin{equation}
    \mathbf{C}_M^{a,\mathrm{fLGR}}
    \begin{bmatrix}
        \mathbf{V}_k\\
        \mathbf{X}_k^{c}\\
        X_{a,k}
    \end{bmatrix}
    =
    \mathbf{0}
    \qquad
    \mathbf{C}_M^{a,\mathrm{fLGR}}
    =
    \begin{bmatrix}
        \mathbf{B}_M^{a,\mathrm{fLGR}} & -\mathbf{I}_M & \mathds{1}_M
    \end{bmatrix}
    \label{eq:flgr-birkhoff-reconstruction-system}
\end{equation}
where the operator's linearity guarantees that its Jacobian is the constant $\mathbf{C}_M^{a,\mathrm{fLGR}}$.

\subsubsection{Connection with the Lagrange Formulation}
\label{subsubsec:flgr-lagrange-connection}

\par
\begin{figure*}[htb]
    \centering
    \begin{subfigure}{3.25in}
        \includegraphics[width=\linewidth]{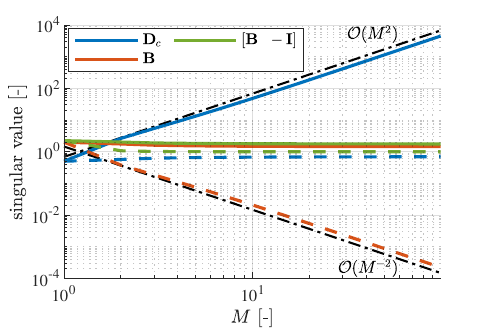}
        \caption{}
        \label{subfig:flgr-singular-values}
    \end{subfigure}%
    \hfill
    \begin{subfigure}{3.25in}
        \includegraphics[width=\linewidth]{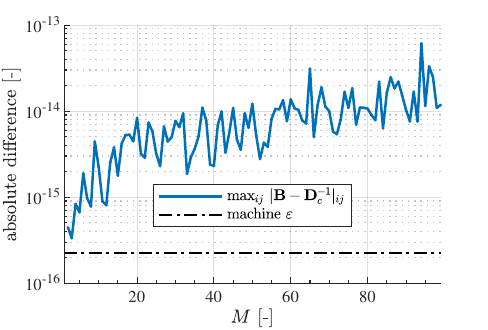}
        \caption{}
        \label{subfig:flgr-birkhoff-lagrange-difference}
    \end{subfigure}
    \caption{Lagrange-Birkhoff connection on \ac{flgr} grids: a) singular values of $\mathbf{D}_c$, $\mathbf{B}_M^{a,\mathrm{fLGR}}$, and $[\mathbf{B}_M^{a,\mathrm{fLGR}}\ \ -\mathbf{I}_M]$; b) maximum deviation of $\mathbf{B}_M^{a,\mathrm{fLGR}}$ from $\mathbf{D}_c^{-1}$.}
    \label{fig:flgr-lagrange-birkhoff-spectrum}
\end{figure*}

 The explicit operator of Eq.~\eqref{eq:flgr-birkhoff-matrix} is also directly related to the Lagrange differentiation matrix on the same grid. Partitioning the \ac{flgr} differentiation matrix of Eq.~\eqref{eq:lagrange-differentiation} as
\begin{equation}
    \mathbf{D}=\begin{bmatrix}\mathbf{d}_a & \mathbf{D}_c\end{bmatrix}
    \qquad
    \mathbf{d}_a\in\mathbb{R}^{M},\quad\mathbf{D}_c\in\mathbb{R}^{M\times M}
    \label{eq:flgr-differentiation-partition}
\end{equation}
where $\mathbf{d}_a$ is the column associated with the noncollocated node $\tau=-1$, and $\mathbf{D}_c$ collects the columns associated with the collocation nodes, two properties follow. Since the derivative of a constant is null, $\mathbf{D}\mathds{1}_{M+1}=\mathbf{0} \rightarrow \mathbf{d}_a=-\mathbf{D}_c\mathds{1}_M$. Moreover, the matrix $\mathbf{D}_c$ is nonsingular, since a polynomial of degree $M$ that vanishes at $\tau=-1$ and whose derivative is null at the $M$ collocation nodes is identically zero. The Lagrange defects of Eq.~\eqref{eq:lagrange-defects} can therefore be rewritten in the integral form $\mathbf{X}_k^c=X_{a,k}\mathds{1}_M+\mathbf{D}_c^{-1}\mathbf{V}_k$, where $\mathbf{D}_c^{-1}$ acts as the implicit integration matrix of the Lagrange \ac{flgr} formulation. Since both formulations describe the unique polynomial of degree $M$ with value $X_{a,k}$ at $\tau=-1$ and derivative values $\mathbf{V}_k$ at the collocation nodes, it holds that
\begin{equation}
    \mathbf{B}_M^{a,\mathrm{fLGR}}=\mathbf{D}_c^{-1}
    \label{eq:flgr-birkhoff-lagrange-equivalence}
\end{equation}
which is the \ac{flgr} counterpart of the equivalence between the differential and integral forms of Radau collocation \cite{garg2010_UnifiedFrameworkNumerical}. Figure~\ref{fig:flgr-lagrange-birkhoff-spectrum} summarizes the numerical verification of this connection for $M=2,\ldots,100$, which prevents the figure scale from masking away relevant features (the same features are observed up to $M=10000$). Figure~\ref{subfig:flgr-birkhoff-lagrange-difference} compares $\mathbf{B}_M^{a,\mathrm{fLGR}}$, computed through Eq.~\eqref{eq:flgr-birkhoff-matrix}, with the numerically inverted $\mathbf{D}_c$, showing a maximum absolute difference in the order of machine precision (dash-dotted line) over the whole range of polynomial degrees. Consequently, for identical grids, the Lagrange and Birkhoff \ac{flgr} transcriptions define the same discrete dynamics and differ only in the structure and conditioning of the resulting \ac{nlp}, while Eq.~\eqref{eq:flgr-birkhoff-matrix} provides the operator explicitly, without matrix inversion.

\par The equivalence in Eq.~\eqref{eq:flgr-birkhoff-lagrange-equivalence} also provides an interpretation of the different conditioning of Lagrange and Birkhoff transcriptions recalled in Sec.~\ref{subsec:birkhoff-overview} and observed in Sec.~\ref{subsec:orbit-raising-problem}. Since $\mathbf{B}_M^{a,\mathrm{fLGR}}=\mathbf{D}_c^{-1}$, the singular values of the two operators are reciprocal, $\sigma_i(\mathbf{B}_M^{a,\mathrm{fLGR}})=1/\sigma_{M+1-i}(\mathbf{D}_c)$, and thus the two matrices share the same condition number. Their ill-conditioning, however, originates at opposite ends of the spectrum, as shown in Fig.~\ref{subfig:flgr-singular-values} (solid and dashed lines denote the largest and smallest singular values, respectively, and black lines the $\mathcal{O}(M^{2})$ and $\mathcal{O}(M^{-2})$ references): the growth of $\kappa_2(\mathbf{D}_c)$ is driven by its largest singular value, which follows the $\mathcal{O}(M^2)$ trend, whereas the largest singular value of $\mathbf{B}_M^{a,\mathrm{fLGR}}$ remains bounded and its smallest one decreases as $\mathcal{O}(M^{-2})$. This distinction is what the lifted Birkhoff representation exploits through the virtual derivative variables (i.e., the identity block in Eq.~\eqref{eq:flgr-birkhoff-reconstruction-system}): appending columns to a matrix adds a positive semidefinite term to its Gram matrix, and since $[\mathbf{B}\ \ -\mathbf{I}_M][\mathbf{B}\ \ -\mathbf{I}_M]^T=\mathbf{B}\mathbf{B}^T+\mathbf{I}_M$, the eigenvalues of this Gram matrix (which correspond to the squared singular values of the original $[\mathbf{B}\ \ -\mathbf{I}_M]$ matrix) cannot decrease below 1 \cite{horn2013_MatrixAnalysis}. The identity block therefore lifts the vanishing singular values of $\mathbf{B}_M^{a,\mathrm{fLGR}}$ away from zero while keeping the largest ones bounded, as shown in Fig.~\ref{subfig:flgr-singular-values} for the endpoint-reduced system $[\mathbf{B}_M^{a,\mathrm{fLGR}}\ \ -\mathbf{I}_M]$. The same mechanism cannot reduce the growing largest singular values of the differentiation matrix, since augmenting $\mathbf{D}_c$ would only increase them further \cite{horn2013_MatrixAnalysis}: this provides an algebraic foundation for the superior conditioning of Birkhoff integration-based formulations with respect to the Lagrange differentiation-based counterparts.

\subsubsection{Reflection from Standard-Radau Operators}
\label{subsubsec:flgr-birkhoff-reflection}

\par It is known that flipped-Radau quadrature weights can be recovered from their direct-Radau counterparts by reflecting their pseudotime domain \cite{sagliano2017_RadauPseudospectralMethod}. Since the Birkhoff \ac{ps} operators are built by integrating the same cardinal polynomials, an analogous reflection relation holds for the Birkhoff integration matrices, as shown below. Introducing the reversal matrix $\mathbf{J}_M\in\mathbb{R}^{M\times M}$, with entries $(\mathbf{J}_M)_{ij}=\delta_{i,M+1-j}$, the \ac{lgr} and \ac{flgr} grids satisfy $\boldsymbol{\tau}^{\mathrm{fLGR}}=-\mathbf{J}_M\boldsymbol{\tau}^{\mathrm{LGR}}$, and the associated cardinal polynomials
\begin{equation}
    \ell_j^{\mathrm{fLGR}}(\tau)
    =
    \ell_{M+1-j}^{\mathrm{LGR}}(-\tau)
    \label{eq:flgr-lgr-cardinal-reflection}
\end{equation}

\par Applying the reflection operation in the context of a Birkhoff formulation also implies changing the endpoint-anchor. In particular, the Birkhoff-$b$ version (i.e., anchored on the final value $\mathbf{X}_b\approx\mathbf{x}(1)$) of the basis functions is given by  $B_j^{b,\mathrm{LGR}}(\tau)=\int_{1}^{\tau}\ell_j^{\mathrm{LGR}}(s)\,\mathrm{d}s$. Notice that the integration is now "oriented" from the final point $\tau=1$ backward. The definition of the Birkhoff-$b$ basis for the \ac{lgr} grid can be recovered by a change of integration variable in Eq.~\eqref{eq:birkhoff-basis-integral}
\begin{equation}
    B_j^{a,\mathrm{fLGR}}(\tau)
    =
    -B_{M+1-j}^{b,\mathrm{LGR}}(-\tau)
    \label{eq:flgr-lgr-birkhoff-basis-reflection}
\end{equation}
and, after evaluation at the corresponding nodes,
\begin{equation}
    \mathbf{B}^{a,\mathrm{fLGR}}
    =
    -\mathbf{J}_M
    \mathbf{B}^{b,\mathrm{LGR}}
    \mathbf{J}_M
    \label{eq:flgr-lgr-birkhoff-matrix-reflection}
\end{equation}
Thus, the \ac{flgr} Birkhoff-$a$ integration matrix can be built directly via Eq.~\eqref{eq:flgr-birkhoff-matrix}, or indirectly by signed reversal of the \ac{lgr} Birkhoff-$b$ matrix included in the universal Birkhoff theory of Refs.~\cite{ross2023_UniversalBirkhoffPseudospectral, ross2024_UniversalBirkhoffTheory}. 
\par The reflection also determines the conditioning of the \ac{flgr} operator. Since $\mathbf{J}_M$ is orthogonal and $\mathbf{J}_M\mathds{1}_M=\mathds{1}_M$, Eq.~\eqref{eq:flgr-lgr-birkhoff-matrix-reflection} gives
\begin{equation}
    \mathbf{C}_M^{a,\mathrm{fLGR}}\big(\mathbf{C}_M^{a,\mathrm{fLGR}}\big)^{T}
    =
    \mathbf{J}_M\left(\mathbf{B}^{b,\mathrm{LGR}}\big(\mathbf{B}^{b,\mathrm{LGR}}\big)^{T}+\mathbf{I}_M+\mathds{1}_M\mathds{1}_M^{T}\right)\mathbf{J}_M
    \label{eq:flgr-lgr-gram-reflection}
\end{equation}
where the term in parentheses is the Gram matrix of the \ac{lgr} Birkhoff-$b$ reconstruction system $[\mathbf{B}^{b,\mathrm{LGR}}\ \ -\mathbf{I}_M\ \ \mathds{1}_M]$. The two systems therefore share their singular values, and the same holds for their endpoint-reduced versions, in which the $\mathds{1}_M$ column is removed. The \ac{flgr} Birkhoff-$a$ operator thus inherits the conditioning recalled in Sec.~\ref{subsec:birkhoff-overview} for Birkhoff transcriptions: $\kappa_2(\mathbf{C}_M^{a,\mathrm{fLGR}})=\sigma_{\max}/\sigma_{\min}$ grows as $\mathcal{O}(\sqrt{M})$ when the anchor value is an unknown, and remains $\mathcal{O}(1)$ when the prescribed endpoint value is eliminated, as verified numerically in Sec.~\ref{subsec:orbit-raising-problem}. This motivates the choice of the Birkhoff-$a$ formulation for the \ac{flgr} grid: since most aerospace \acp{ocp} have known initial conditions, the endpoint-reduced version of Eq.~\eqref{eq:flgr-birkhoff-reconstruction-system} is obtained by simply moving the $X_{a,k}$ term to the known right-hand side. The conditioning of the complete \ac{nlp} and of its associated \ac{kkt} system, however, also depends on the nonlinear problem-specific functions and on the underlying optimization algorithm, an aspect addressed in this work through the scaling approach of Sec.~\ref{subsec:composable-scaling}.

%% file: sections/differential_constraints.tex
\section{Differential Constraints}
\label{sec:differential-constraints}

\par Differential constraints represent bounds on the time derivative of a quantity of interest. They arise when limiting actuator or guidance-command rates, load variations, heating-rate changes, and related quantities, and can also reformulate nonlinear dynamics in control-affine form \cite{wang2017_ConstrainedTrajectoryOptimization}. Consider a continuously differentiable quantity
\begin{equation}
    \mathbf{y}(t)=\mathbf{h}\bigl(\mathbf{x}(t),\mathbf{u}(t),t\bigr)
    \label{eq:differential-quantity-definition}
\end{equation}
where $\mathbf{h}$ may be any nonlinear function. Its physical-time derivative is constrained by
\begin{equation}
    \dot{\mathbf{y}}_{L}(t)\leq\frac{\mathrm{d}\mathbf{y}}{\mathrm{d}t}\leq\dot{\mathbf{y}}_{U}(t)
    \label{eq:continuous-differential-constraints}
\end{equation}
Multiple prior works handle this class of constraints through state augmentation \cite{wang2017_ConstrainedTrajectoryOptimization,wang2020_ImprovedSequentialConvex}. Differentiation-based \ac{ps} formulations also provide a straightforward way to enforce them directly \cite{nie2021_EfficientImplementationRate}, but the same cannot be said for integration-based methods, where derivative quantities are not linear functions of the nodal values and require the introduction of auxiliary variables (e.g., Eq.~\eqref{eq:birkhoff-derivative-variables}). This section extends the direct enforcement of differential constraints to Birkhoff transcriptions, and establishes both the primal and the dual equivalence of the direct approach with state augmentation.

\par Under the time mapping defined in Eq.~\eqref{eq:time-dilation-definition}, with positive time-dilation factor $\rho(\tau)=\mathrm{d}t/\mathrm{d}\tau>0$, we may re-arrange
\begin{equation}
    \frac{\mathrm{d}\mathbf{y}}{\mathrm{d}\tau}=\rho(\tau)\frac{\mathrm{d}\mathbf{y}}{\mathrm{d}t} \quad \Rightarrow \quad \dot{\mathbf{y}}_{L}(t)\leq\frac{1}{\rho(\tau)}\frac{\mathrm{d}\mathbf{y}}{\mathrm{d}\tau}\leq\dot{\mathbf{y}}_{U}(t)
    \label{eq:pseudotime-differential-bounds}
\end{equation}
As done for the state trajectory in Eq.~\eqref{eq:lagrange-state-interpolation}, we approximate the constrained quantity $\mathbf{y}$ from its values at the $N_d$ discretization nodes through the Lagrange interpolant
\begin{equation}
    \mathbf{y}^{L}(\tau)=\sum_{i=1}^{N_d}\mathbf{Y}_i\ell_i(\tau)
    \qquad
    \mathbf{Y}_i\approx\mathbf{y}\bigl(t(\tau_i^d)\bigr)
    =\mathbf{h}\bigl(\mathbf{x}(\tau_i^d),\mathbf{u}(\tau_i^d),t(\tau_i^d)\bigr)
    \label{eq:lagrange-differential-quantity}
\end{equation}
This representation introduces an assumption that the standard transcription does not require. A \ac{ps} method defines the control only at the collocation nodes and makes no assumption on its continuous-time trajectory \cite{ross2015_PrimerPontryaginsPrinciple}. Approximating $\mathbf{y}$ through an interpolant, instead, assigns to $\mathbf{y}$ (and, through $\mathbf{h}$, implicitly to the controls) the same polynomial representation adopted for the states. At the noncollocated node $\tau=-1$ of the \ac{flgr} grid, $\mathbf{h}$ is evaluated using the additional control $\mathbf{U}_a$ introduced in Sec.~\ref{subsec:new-birkhoff-interpolants}, which is then determined by the optimizer rather than extrapolated. This mirrors the free initial value of the augmented state in state augmentation, and enables the primal equivalence of Sec.~\ref{subsec:control-rate-primal-equivalence}.
Using the differentiation operator defined in Eq.~\eqref{eq:lagrange-differentiation}, the differential bounds in Eq.~\eqref{eq:pseudotime-differential-bounds} are imposed as
\begin{equation}
    \dot{\mathbf{y}}_{L}(t_j)
    \leq\frac{1}{\rho_j}\sum_{i=1}^{N_d}D_{ji}\mathbf{Y}_i
    \leq\dot{\mathbf{y}}_{U}(t_j)
    \qquad j=1,\ldots,N_c
    \label{eq:lagrange-differential-constraints}
\end{equation}

\par A Birkhoff-$a$ formulation for the same constraint is also available. Introducing the virtual derivative values $\mathbf{W}_j\approx\mathrm{d}\mathbf{y}/\mathrm{d}\tau|_{\tau_j^c}$ and using the integration matrix defined in Eq.~\eqref{eq:birkhoff-reconstruction}, we impose
\begin{equation}
    \mathbf{Y}_i^c=\mathbf{Y}_a+\sum_{j=1}^{N_c}B_{ij}\mathbf{W}_j
    \qquad i=1,\ldots,N_c
    \label{eq:birkhoff-differential-reconstruction}
\end{equation}
where $\mathbf{Y}_i^c$ is obtained by evaluating $\mathbf{h}$ at the $i$-th collocation node, and $\mathbf{Y}_a=\mathbf{h}(\mathbf{X}_a,\mathbf{U}_a,t_0)$ is evaluated at the noncollocated initial node, playing for $\mathbf{y}$ the same role as $\mathbf{X}_a$ for the states, so that the representation assumption discussed above also applies to the Birkhoff interpolant of Eq.~\eqref{eq:birkhoff-interpolant}.
The derivative bounds are then imposed directly on the virtual derivative values as
\begin{equation}
    \dot{\mathbf{y}}_{L}(t_j)
    \leq\frac{1}{\rho_j}\mathbf{W}_j
    \leq\dot{\mathbf{y}}_{U}(t_j)
    \qquad j=1,\ldots,N_c
    \label{eq:birkhoff-differential-bounds}
\end{equation}
Both formulations only require the evaluation of $\mathbf{h}$ at the discretization nodes, without explicitly defining $\dot{\mathbf{y}}$. Control-rate constraints correspond to the particular choice $\mathbf{h}(\mathbf{x},\mathbf{u},t)=\mathbf{u}$, so that $\dot{\mathbf{y}}=\dot{\mathbf{u}}$: since these are common in aerospace \acp{ocp}, the following subsections compare their direct enforcement with state augmentation \cite{wang2017_ConstrainedTrajectoryOptimization}.

\subsection{Primal Equivalence to State Augmentation}
\label{subsec:control-rate-primal-equivalence}

\par We start by analyzing the relations between primal variables defined by the two techniques. For differentiation-based formulations, this result is consistent with Ref.~\cite{nie2021_EfficientImplementationRate} and is reproduced here for completeness, as it underlies both its Birkhoff counterpart and the dual analysis of Sec.~\ref{subsec:control-rate-dual-equivalence}. Consider a control trajectory $\mathbf{u}(t)\in\mathbb{R}^{n_u}$ subject to the physical-time rate bounds
\begin{equation}
    \underline{\mathbf{v}}(t)\leq\dot{\mathbf{u}}(t)\leq\overline{\mathbf{v}}(t)
    \label{eq:continuous-control-rate-limits}
\end{equation}
where $\underline{\mathbf{v}}$ and $\overline{\mathbf{v}}$ denote the lower and upper admissible control rates. For nodal control values $\mathbf{U}_i\approx\mathbf{u}(t(\tau_i^d))$, the direct Lagrange formulation in Eq.~\eqref{eq:lagrange-differential-constraints} becomes
\begin{equation}
    \underline{\mathbf{v}}(t_j)
    \leq\frac{1}{\rho_j}\sum_{i=1}^{N_d}D_{ji}\mathbf{U}_i
    \leq\overline{\mathbf{v}}(t_j)
    \qquad j=1,\ldots,N_c
    \label{eq:direct-control-rate-transcription}
\end{equation}
The corresponding Birkhoff form follows directly from Eqs.~\eqref{eq:birkhoff-differential-reconstruction} and~\eqref{eq:birkhoff-differential-bounds} with $\mathbf{Y}_i^c=\mathbf{U}_i$ and $\mathbf{W}_j\approx\mathrm{d}\mathbf{u}/\mathrm{d}\tau|_{\tau_j^c}$.

\par The state-augmentation formulation instead promotes $\mathbf{u}$ to a state and introduces its physical-time derivative $\mathbf{v}$ as a new control
\begin{equation}
    \begin{aligned}
        \dot{\mathbf{u}}(t)=\mathbf{v}(t) &\qquad
        \underline{\mathbf{v}}(t)\leq\mathbf{v}(t)\leq\overline{\mathbf{v}}(t)
    \end{aligned}
    \label{eq:control-rate-state-augmentation}
\end{equation}
Let $\mathbf{V}_j\approx\mathbf{v}(t_j)$ denote the nodal values of this additional control. Using the same grid, interpolation, and time mapping as the direct formulation, the added dynamics are transcribed as
\begin{equation}
    \sum_{i=1}^{N_d}D_{ji}\mathbf{U}_i=\rho_j\mathbf{V}_j
    \qquad j=1,\ldots,N_c
    \label{eq:augmented-control-dynamics-defect}
\end{equation}
Since $\rho_j>0$, each collocation equation can be solved for the rate control
\begin{equation}
    \mathbf{V}_j=\frac{1}{\rho_j}\sum_{i=1}^{N_d}D_{ji}\mathbf{U}_i
    \qquad j=1,\ldots,N_c
    \label{eq:eliminated-rate-control}
\end{equation}
Substitution into the bounds on $\mathbf{V}_j$ recovers Eq.~\eqref{eq:direct-control-rate-transcription}. Therefore, for identical grids, interpolation, endpoint treatment (with $\mathbf{U}_a$ corresponding to the initial value of the augmented state), and time mapping, both formulations define the same feasible set for the nodal controls and consequently admit the same primal solutions \cite{nie2021_EfficientImplementationRate}. The direct formulation is the algebraically condensed form of state augmentation: it eliminates the nodal rate controls $\mathbf{V}_j$ and their associated collocation equations before the optimization problem is solved. It therefore represents the same constraints with fewer decision variables and equality constraints, providing a more compact formulation and a potential computational advantage (note that the computational performance also depends on the problem structure and underlying \ac{nlp} solver). The same result holds for the Birkhoff formulation, where $\mathbf{W}_j/\rho_j$ plays the role of the eliminated rate control $\mathbf{V}_j$.

\subsection{Dual Equivalence to State Augmentation}
\label{subsec:control-rate-dual-equivalence}

\par A second interpretation of this equivalence can be recovered through the first-order optimality conditions. If the original running cost does not depend explicitly on $\dot{\mathbf{u}}$, state augmentation treats $\mathbf{u}$ as a state and $\mathbf{v}=\dot{\mathbf{u}}$ as its control \cite{wang2017_ConstrainedTrajectoryOptimization, wang2020_ImprovedSequentialConvex}. The Hamiltonian of the state-augmented \ac{ocp} is
\begin{equation}
    \mathcal{H}_{A}
    =
    L(\mathbf{x},\mathbf{u},t)
    +
    \boldsymbol{\lambda}_{x}^{T}\mathbf{f}(\mathbf{x},\mathbf{u},t)
    +
    \boldsymbol{\lambda}_{u}^{T}\mathbf{v}
    \label{eq:augmented-control-rate-hamiltonian}
\end{equation}
where $\boldsymbol{\lambda}_{x}$ and $\boldsymbol{\lambda}_{u}$ are the costates associated with $\mathbf{x}$ and with the augmented state $\mathbf{u}$, respectively, and $L(\mathbf{x},\mathbf{u},t)$ is the running cost. Since $\mathbf{v}$ enters Eq.~\eqref{eq:augmented-control-rate-hamiltonian} linearly, the augmentation introduces an affine dependence on the auxiliary control, even when the original control $\mathbf{u}$ enters the dynamics regularly. The corresponding switching function is $\partial\mathcal{H}_{A}/\partial\mathbf{v}=\boldsymbol{\lambda}_{u}$, so the optimal rate is generally driven to one of its admissible bounds, while intervals on which the relevant component of $\boldsymbol{\lambda}_{u}$ vanishes are singular arcs and require higher-order optimality conditions to determine the rate \cite{ross2015_PrimerPontryaginsPrinciple, longuski2014_OptimalControlAerospace}.

\par The adjoint equation associated with the augmented state $\mathbf{u}$ follows directly from the Pontryagin costate equation
\begin{equation}
    \dot{\boldsymbol{\lambda}}_{u}
    =
    -\frac{\partial\mathcal{H}_{A}}{\partial\mathbf{u}}
    =
    -\left(
        \frac{\partial L}{\partial\mathbf{u}}
        +
        \frac{\partial\mathbf{f}^{T}}{\partial\mathbf{u}}
        \boldsymbol{\lambda}_{x}
    \right)
    \label{eq:control-rate-auxiliary-costate}
\end{equation}
where the term $\boldsymbol{\lambda}_{u}^{T}\mathbf{v}$ from Eq.~\eqref{eq:augmented-control-rate-hamiltonian} does not contribute because $\mathbf{v}$ is an independent control in the augmented formulation. Together with the minimum condition on $\mathbf{v}$, Eq.~\eqref{eq:control-rate-auxiliary-costate} determines the switching structure of the rate-limited solution.

\par The same adjoint structure can be recovered when $\mathbf{u}$ remains a control subject to direct differential constraints. We introduce the nonnegative multipliers $\boldsymbol{\mu}_{L}$ and $\boldsymbol{\mu}_{U}$ for the inequalities $\underline{\mathbf{v}}-\dot{\mathbf{u}}\leq\mathbf{0}$ and $\dot{\mathbf{u}}-\overline{\mathbf{v}}\leq\mathbf{0}$, respectively. Adjoining these differential constraints to the Hamiltonian gives the augmented Hamiltonian 
\begin{equation} \overline{\mathcal{H}}_{D} = L(\mathbf{x},\mathbf{u},t) + \boldsymbol{\lambda}_{x}^{T}\mathbf{f}(\mathbf{x},\mathbf{u},t) + \boldsymbol{\mu}_{L}^{T}\left(\underline{\mathbf{v}}-\dot{\mathbf{u}}\right) + \boldsymbol{\mu}_{U}^{T}\left(\dot{\mathbf{u}}-\overline{\mathbf{v}}\right)
\label{eq:direct-rate-augmented-hamiltonian} 
\end{equation} 
Unlike a standard path-constrained Hamiltonian, $\overline{\mathcal{H}}_{D}$ depends explicitly on $\dot{\mathbf{u}}$. The stationarity condition $\frac{\mathrm{d}\overline{\mathcal{H}}_{D}}{\mathrm{d}\mathbf{u}}=0$ therefore requires the expansion of the total derivative into
\begin{equation} 
\frac{\partial\overline{\mathcal{H}}_{D}}{\partial\mathbf{u}} - \frac{\mathrm{d}}{\mathrm{d}t} \left( \frac{\partial\overline{\mathcal{H}}_{D}}{\partial\dot{\mathbf{u}}} \right) = \mathbf{0}
\label{eq:direct-rate-generalized-stationarity} 
\end{equation} 
We now compute each term in Eq.~\eqref{eq:direct-rate-generalized-stationarity}, starting from the derivative of Eq.~\eqref{eq:direct-rate-augmented-hamiltonian} with respect to $\dot{\mathbf{u}}$:
\begin{equation}
\frac{\partial\overline{\mathcal{H}}_{D}}{\partial\dot{\mathbf{u}}} = -\boldsymbol{\mu}_{L} + \boldsymbol{\mu}_{U}
\label{eq:direct-rate-derivative-multiplier}
\end{equation}
We now define the variable $\boldsymbol{\lambda}_{u}$ as
\begin{equation} \boldsymbol{\lambda}_{u} \triangleq \boldsymbol{\mu}_{L} - \boldsymbol{\mu}_{U}
\label{eq:implicit-control-costate} 
\end{equation} 
Note that no implication that this $\boldsymbol{\lambda}_{u}$ corresponds to the one for the state-augmented system in Eq.~\eqref{eq:augmented-control-rate-hamiltonian} has been made yet. Substituting this definition into Eq.~\eqref{eq:direct-rate-derivative-multiplier} yields
\begin{equation}
\frac{\partial\overline{\mathcal{H}}_{D}}{\partial\dot{\mathbf{u}}} = -\boldsymbol{\lambda}_{u}
\label{eq:direct-rate-implicit-costate-derivative}
\end{equation} 
The missing term in Eq.~\eqref{eq:direct-rate-generalized-stationarity} is
\begin{equation}
\frac{\partial\overline{\mathcal{H}}_{D}}{\partial\mathbf{u}} = \frac{\partial L}{\partial\mathbf{u}} + \frac{\partial\mathbf{f}^{T}}{\partial\mathbf{u}} \boldsymbol{\lambda}_{x}
\label{eq:direct-rate-control-gradient}
\end{equation}
Substitution of the known terms of Eq.~\eqref{eq:direct-rate-generalized-stationarity} therefore yields 
\begin{equation}
\dot{\boldsymbol{\lambda}}_{u} = -\left( \frac{\partial L}{\partial\mathbf{u}} + \frac{\partial\mathbf{f}^{T}}{\partial\mathbf{u}} \boldsymbol{\lambda}_{x} \right) 
\label{eq:direct-rate-implicit-costate}
\end{equation} 
which is identical to the auxiliary costate equation obtained from the state-augmentation formulation in Eq.~\eqref{eq:control-rate-auxiliary-costate}. The quantity $\boldsymbol{\lambda}_{u}=\boldsymbol{\mu}_{L}-\boldsymbol{\mu}_{U}$ therefore plays in the direct formulation the same adjoint role as the explicit costate associated with the augmented state $\mathbf{u}$. The primal and dual equivalence results show that the direct formulation preserves the same auxiliary adjoint structure and switching structure as the state augmentation approach, while eliminating explicit rate control from the primal \ac{ocp}.

%% file: sections/numerical_advancements.tex
\section{Numerical Advancements}
\label{sec:numerical-advancements}

\par The theoretical developments in Secs.~\ref{sec:birkhoff-theory} and~\ref{sec:differential-constraints} are integrated into a common numerical framework, implemented in \ac{raptor}, that enables heterogeneous multiphase \ac{ps} transcription: each phase is transcribed according to its independently defined node family, polynomial order, and underlying global interpolant (Lagrange or Birkhoff), while endpoint and linkage conditions couple the phasewise discretizations into a single \ac{nlp}.

\subsection{Heterogeneous Multiphase Transcription}
\label{subsec:heterogeneous-multiphase-transcription}

\par Consider an \ac{ocp} composed of $N_{\mathrm{ph}}$ phases indexed by $q=1,\ldots,N_{\mathrm{ph}}$. Phase $q$ is defined on the physical-time interval $[t_0^{(q)},t_f^{(q)}]$ by the state $\mathbf{x}^{(q)}(t)\in\mathbb{R}^{n_x^{(q)}}$, the control $\mathbf{u}^{(q)}(t)\in\mathbb{R}^{n_u^{(q)}}$, and the static parameters $\mathbf{p}^{(q)}\in\mathbb{R}^{n_p^{(q)}}$. The vector $\mathbf{p}_g\in\mathbb{R}^{n_g}$ collects global parameters (i.e., shared among phases). For compactness, the endpoint data of each phase are collected as
\begin{equation}
    \mathbf{e}^{(q)}
    =
    \begin{bmatrix}
        \mathbf{x}^{(q)}(t_0^{(q)})&
        \mathbf{x}^{(q)}(t_f^{(q)})&
        t_0^{(q)}&
        t_f^{(q)}&
        \mathbf{p}^{(q)}
    \end{bmatrix}^{T}
    \label{eq:phase-endpoint-vector}
\end{equation}
and the arguments of the phase functions are shortened as $(\cdot)^{(q)}=(\mathbf{x}^{(q)}(t),\mathbf{u}^{(q)}(t),t,\mathbf{p}^{(q)},\mathbf{p}_g)$. The multiphase \ac{ocp} is then written in Bolza form as
\begin{equation}
    \begin{aligned}
        \operatorname{minimize}\quad
        &J=\sum_{q=1}^{N_{\mathrm{ph}}}\left[\Phi^{(q)}\!\left(\mathbf{e}^{(q)},\mathbf{p}_g\right)+\int_{t_0^{(q)}}^{t_f^{(q)}}L^{(q)}(\cdot)^{(q)}\,\mathrm{d}t\right]\\
        \operatorname{subject\ to}\quad
        &\dot{\mathbf{x}}^{(q)}(t)=\mathbf{f}^{(q)}(\cdot)^{(q)}
        \qquad
        \mathbf{g}_{L}^{(q)}(t)\leq\mathbf{g}^{(q)}(\cdot)^{(q)}\leq\mathbf{g}_{U}^{(q)}(t)\\
        &\mathbf{b}_{L}^{(q)}\leq\mathbf{b}^{(q)}\!\left(\mathbf{e}^{(q)},\mathbf{p}_g\right)\leq\mathbf{b}_{U}^{(q)}
    \end{aligned}
    \label{eq:multiphase-ocp}
\end{equation}
where $\Phi^{(q)}$ and $L^{(q)}$ are the endpoint and running costs of phase $q$, respectively. The dynamics and path constraints hold for the full phase interval $t\in[t_0^{(q)},t_f^{(q)}]$, whereas the endpoint constraints act on the endpoint data, and equality constraints are recovered by defining identical lower and upper bounds. Box constraints can be defined directly on the states, controls, phase times, and static parameters. Differential constraints follow from Eqs.~\eqref{eq:differential-quantity-definition} and~\eqref{eq:continuous-differential-constraints} by applying $\mathbf{h}^{(q)}$ to the variables and parameters of phase $q$.

\par The phases are coupled by nonlinear linkage constraints acting on their endpoint data. Each linkage $r=1,\ldots,N_{\mathrm{link}}$ involves two arbitrary phases $q_1^{(r)}$ and $q_2^{(r)}$, and the corresponding constraint is
\begin{equation}
    \boldsymbol{\eta}_{L}^{(r)}
    \leq
    \boldsymbol{\eta}^{(r)}\!\left(
        \mathbf{e}^{(q_1^{(r)})},
        \mathbf{e}^{(q_2^{(r)})},
        \mathbf{p}_g
    \right)
    \leq
    \boldsymbol{\eta}_{U}^{(r)}
    \label{eq:multiphase-linkage-constraints}
\end{equation}
Since $q_1$ and $q_2$ need not be consecutive phase indexes and no assumption is made on the number of linkages $N_{\mathrm{link}}$, a linkage is not restricted to only enforce continuity conditions between subsequent phases. These features open several opportunities for modeling complex systems, including state and time continuity even through nonlinear transformations such as coordinate changes \cite{sagliano2025_SixDegreeofFreedomEntryPowered}, enforcing discontinuities between phases, and branching trajectories by linking the terminal conditions of a common phase with the initial conditions of multiple successor phases, as demonstrated in the \ac{rlv} problem in Sec.~\ref{subsec:rlv-problem}. 

\subsection{Variable-Dependent Non-Affine Time Mappings}
\label{subsec:non-affine-time-mappings}

\par The heterogeneous treatment of each \ac{ocp} phase enables targeted adjustments of both the \ac{ps} operators and the associated time mapping. Omitting the phase index for clarity, the time mapping transformation $\mathcal{T}$ of Eq.~\eqref{eq:time-dilation-definition} can also be expressed through a normalized mapping function $s$ as
\begin{equation}
    \begin{aligned}
        t
        &=
        \mathcal{T}(\tau;t_0,t_f,\mathbf{p})
        =
        t_0+\left(t_f-t_0\right)s\!\left(\tau;\mathbf{p}\right)
        \\
        \rho(\tau;t_0,t_f,\mathbf{p})
        &=
        \frac{\partial\mathcal{T}}{\partial\tau}(\tau;t_0,t_f,\mathbf{p})
        =
        \left(t_f-t_0\right)
        \frac{\partial s}{\partial\tau}\!\left(\tau;\mathbf{p}\right)
    \end{aligned}
    \label{eq:normalized-time-mapping}
\end{equation}
where $s$ defines how each pseudotime value $\tau$ is mapped onto the normalized interval $[0, 1]$, and Eq.~\eqref{eq:normalized-time-mapping} then scales it onto the physical-time interval $[t_0, t_f]$. Required conditions for a valid non-affine time mapping are $s(-1;\mathbf{p})=0$, $s(1;\mathbf{p})=1$, and $\partial s/\partial\tau>0$: these preserve the phase endpoints and ensure a monotonic transformation. The conventional affine map is recovered with $s=(1+\tau)/2$, which produces the constant dilation $\rho=(t_f-t_0)/2$.

\par A non-affine choice of $s$ changes the distribution of the fixed \ac{ps} nodes in physical time without modifying the underlying node family or interpolation operators \cite{ross2007_LowThrustHighAccuracyTrajectory, ross2020_EnhancementsDIDOOptimal}. Regions of smaller dilation contain more physical-time nodes and therefore receive greater local resolution, whereas regions of larger dilation are represented more coarsely. This provides a direct means of resolving localized trajectory features, such as rapid transients or control switches, while retaining the characteristic global polynomial approximation of \ac{ps} methods. An example of this approach is shown in Sec.~\ref{subsubsec:landing-time-mapping}.

\par The time mapping must be applied consistently throughout the transcription. The state dynamics are transformed according to Eq.~\eqref{eq:pseudotime-dynamics}, and the corresponding transformation of differential constraints is given by Eq.~\eqref{eq:pseudotime-differential-bounds}. The running-cost contribution becomes
\begin{equation}
    \int_{t_0}^{t_f}
    L\!\left(\mathbf{x}(t),\mathbf{u}(t),t,\mathbf{p},\mathbf{p}_g\right)
    \,\mathrm{d}t
    =
    \int_{-1}^{1}
    \rho(\tau;t_0,t_f,\mathbf{p})
    L\!\left(\mathbf{x}(\tau),\mathbf{u}(\tau),\mathcal{T}(\tau;t_0,t_f,\mathbf{p}),\mathbf{p},\mathbf{p}_g\right)
    \,\mathrm{d}\tau
    \label{eq:time-mapped-running-cost}
\end{equation}
When $t_0$, $t_f$, or components of $\mathbf{p}$ entering $s$ are optimization variables, the physical location of every collocation node and the corresponding dilation vary according to the decision variables. This dependence enters the dynamics both through the factor $\rho$ multiplying $\mathbf{f}$ and through the mapped time $\mathcal{T}$ appearing as an argument of $\mathbf{f}$ (e.g., in non-autonomous systems \cite{harris2014_LosslessConvexificationNonconvex}), and analogously the time-dependent path and differential constraints and the running cost of Eq.~\eqref{eq:time-mapped-running-cost}: retaining these chain-rule terms is therefore necessary for the derivatives of the transcribed problem to remain consistent with the selected time mapping. The time dilation factor is also required for consistent application of the covector mapping theorem associated with the \ac{ps} operators used in each phase \cite{gong2008_ConnectionsCovectorMapping, ross2001_PseudospectralTransformationCovectors}.

\subsection{Composable Scaling and Equilibration}
\label{subsec:composable-scaling}

\par This section reviews scaling methods that improve the conditioning of the transcribed \ac{ocp}, and shows how they can be combined to recover existing scaling methods for \acp{ocp} or to design new ones. We define $\mathbf{z}$ as the \ac{nlp} decision vector, $\mathbf{c}(\mathbf{z})$ the constraint vector, and $J(\mathbf{z})$ the objective. Their scaled counterparts are defined by
\begin{equation}
    \widehat{\mathbf{z}}=\mathbf{A}_z\mathbf{z}+\mathbf{b}_z\qquad \widehat{\mathbf{c}}=\mathbf{A}_c\mathbf{c}\qquad \widehat{J}=s_JJ
    \label{eq:scaling-maps}
\end{equation}
where $\mathbf{A}_z$ and $\mathbf{A}_c$ are nonsingular matrices, $\mathbf{b}_z$ translates the decision variables, and $s_J>0$ scales the objective. Constraint bounds are transformed by the same matrix $\mathbf{A}_c$. If these transformations are nonsingular, the physical feasible set and minimizers are preserved. Using column-vector gradients, the transformed objective gradient and constraint Jacobian are, respectively,
\begin{equation}
    \nabla_{\widehat{\mathbf{z}}}\widehat{J}=s_J\mathbf{A}_z^{-T}\nabla_{\mathbf{z}}J\qquad \mathbf{J}_{\widehat{c}}=\mathbf{A}_c\mathbf{J}_{c}\mathbf{A}_z^{-1}
    \label{eq:scaled-derivatives}
\end{equation}
Analogously, defining $\widehat{\boldsymbol{\lambda}}$ and $\widehat{\boldsymbol{\nu}}$ as the constraint and variable-bound multipliers in scaled coordinates, their physical-coordinate counterparts are
\begin{equation}
    \boldsymbol{\lambda}=\frac{1}{s_J}\mathbf{A}_c^{T}\widehat{\boldsymbol{\lambda}}\qquad \boldsymbol{\nu}=\frac{1}{s_J}\mathbf{A}_z^{T}\widehat{\boldsymbol{\nu}}
    \label{eq:dual-unscaling}
\end{equation}
Equations~\eqref{eq:scaled-derivatives} and~\eqref{eq:dual-unscaling} ensure that derivatives and dual quantities remain consistent when interpreted in physical coordinates.
\par \textit{Manual scaling} assigns characteristic values to individual variables, constraints, and the objective, which directly define the diagonal entries of $\mathbf{A}_z$ and $\mathbf{A}_c$, the translation $\mathbf{b}_z$, and $s_J$ in Eq.~\eqref{eq:scaling-maps}. These characteristic values may be defined with the goal of scaling (i.e., changing the numerical coordinates of the primal problem) or balancing (i.e., targeting compatibility with the associated dual variables and underlying necessary conditions): primal quantities of comparable magnitude do not necessarily yield a balanced problem \cite{ross2018_ScalingBalancingHighPerformance}.

\par \textit{Bounds projection} maps each variable with finite and distinct bounds to the unit interval \cite{betts2010_PracticalMethodsOptimal}. For $z_{L,i}<z_{U,i}$,
\begin{equation}
    (\mathbf{A}_z)_{ii}=\frac{1}{z_{U,i}-z_{L,i}}\qquad (\mathbf{b}_z)_i=-\frac{z_{L,i}}{z_{U,i}-z_{L,i}}\qquad \widehat{z}_i=\frac{z_i-z_{L,i}}{z_{U,i}-z_{L,i}}\in[0,1]
    \label{eq:bounds-projection-scaling}
\end{equation}
Variables that are fixed or without two finite bounds are left unchanged.

\par \textit{\ac{jrn}} scales the constraints and objective according to their derivatives at a reference point, typically the initial guess. Let $\mathbf{G}_0$ and $\mathbf{g}_0$ denote the constraint Jacobian and objective gradient computed on the initial guess $\mathbf{z}_0$. For a selected norm $\|\cdot\|_r$, rows with $0$ norm are left unscaled, while a safeguard $\epsilon_s>0$ is used to compute the additional constraint and objective factors as
\begin{equation}
    (\mathbf{A}_{\mathrm{JRN}})_{ii}=\frac{1}{\max\!\left(\left\|\mathbf{e}_i^{T}\mathbf{G}_0\right\|_r,\epsilon_s\right)}\qquad s_{J,\mathrm{JRN}}=\frac{1}{\max\!\left(\|\mathbf{g}_0\|_r,\epsilon_s\right)}
    \label{eq:jacobian-row-normalization}
\end{equation}
where $\mathbf{e}_i^{T}$ is the $i$-th row selection operator.

\par \textit{Static Ruiz equilibration} iteratively balances the row and column norms of a reference Jacobian $\mathbf{G}^{(0)}$, evaluated at the initial guess after any preceding scaling transformations. Let $\|\cdot\|_r$ denote the selected vector norm. At iteration $k$, rows and columns with 0 norm are left unscaled (i.e., the corresponding entries in the row and column scaling matrices are $(\mathbf{R}^{(k)})_{ii}=(\mathbf{C}^{(k)})_{ii}=1$), while for non-zero row terms the following algorithm applies \cite{ruiz2001_ScalingAlgorithmEquilibrate}
\begin{equation}
        (\mathbf{R}^{(k)})_{ii}
        =
        \!\left(\|\mathbf{e}_i^{T}\mathbf{G}^{(k)}\|_r\right)^{-1/2}
        \qquad
        (\mathbf{C}^{(k)})_{jj}
        =
        \!\left(\|\mathbf{G}^{(k)}\mathbf{e}_j\|_r\right)^{-1/2}
    \label{eq:ruiz-equilibration-factors}
\end{equation}
and the reference Jacobian is updated according to
\begin{equation}
    \mathbf{G}^{(k+1)}=\mathbf{R}^{(k)}\mathbf{G}^{(k)}\mathbf{C}^{(k)}
    \label{eq:ruiz-equilibration-update}
\end{equation}
The accumulated row factors contribute to $\mathbf{A}_c$, while the accumulated column factors contribute to $\mathbf{A}_z^{-1}$. The iteration terminates when the nonzero row and column norms are within a prescribed tolerance of unity, or when the maximum iteration count is reached. The objective scaling factor $s_J$ is then obtained as the row factor of the same iteration applied to the reference objective gradient, transformed by the accumulated column scaling, while the variable scaling obtained from the constraint Jacobian is left unchanged.

\par We now introduce a composable scaling construction, showing how a sequence of affine scaling and equilibration methods can be reduced to a single equivalent transformation. If two variable transformations are applied as $\mathbf{z}^{(1)}=\mathbf{A}_1\mathbf{z}+\mathbf{b}_1$ and $\mathbf{z}^{(2)}=\mathbf{A}_2\mathbf{z}^{(1)}+\mathbf{b}_2$, their composition is defined as
\begin{equation}
    \mathbf{z}^{(2)}=\mathbf{A}_{\mathrm{tot}}\mathbf{z}+\mathbf{b}_{\mathrm{tot}}\qquad \mathbf{A}_{\mathrm{tot}}=\mathbf{A}_2\mathbf{A}_1\qquad \mathbf{b}_{\mathrm{tot}}=\mathbf{A}_2\mathbf{b}_1+\mathbf{b}_2
    \label{eq:composed-variable-scaling}
\end{equation}
while the constraint and objective transformations satisfy
\begin{equation}
    \mathbf{A}_{c,\mathrm{tot}}=\mathbf{A}_{c,2}\mathbf{A}_{c,1}\qquad s_{J,\mathrm{tot}}=s_{J,2}s_{J,1}
    \label{eq:composed-function-scaling}
\end{equation}
The order is significant because derivative-based strategies act on the Jacobian and objective gradient produced by all preceding transformations. We show the application of this approach to the ordered composition of bounds projection, applied through $\mathbf{A}_{\mathrm{BP}}$, and \ac{jrn}, which is computed by replacing $\mathbf{G}_0$ in Eq.~\eqref{eq:jacobian-row-normalization} with the projected Jacobian $\mathbf{G}_{\mathrm{P}}=\mathbf{J}_{c}(\mathbf{z}_0)\mathbf{A}_{\mathrm{BP}}^{-1}$, thus obtaining
\begin{equation}
    \widehat{\mathbf{J}}_c
    =
    \mathbf{A}_{\mathrm{JRN}}\mathbf{G}_{\mathrm{P}}
    =
    \mathbf{A}_{\mathrm{JRN}}\mathbf{J}_{c}(\mathbf{z}_0)\mathbf{A}_{\mathrm{BP}}^{-1}
    \label{eq:pjrn-composition}
\end{equation}
The translation associated with bounds projection does not appear because it does not affect the Jacobian. Note that Eq.~\eqref{eq:pjrn-composition} recovers the \ac{pjrn} strategy introduced in Ref.~\cite{sagliano2017_RadauPseudospectralMethod}. More generally, the composition rule permits variable, constraint, objective, and equilibration strategies to be combined without defining a separate transformation for every sequence, allowing the reconstruction of existing methods as well as the definition of problem-tailored scaling strategies.

\subsection{Hybrid Jacobian Assembly and Solution Recovery}
\label{subsec:hybrid-jacobians}

\par
\begin{figure*}[htb]
    \centering
    \includegraphics[width=\textwidth]{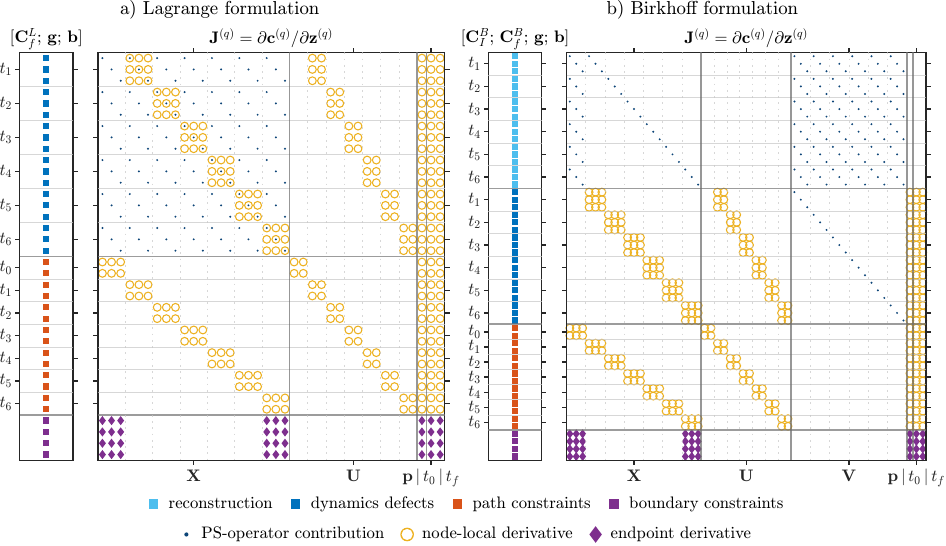}
    \caption{Phase constraint vector and Jacobian structure for a) Lagrange and b) Birkhoff transcriptions.}
    \label{fig:hybrid-jacobian-assembly}
\end{figure*}

 Problem-specific functions such as dynamics, costs, and constraints are typically arbitrary nonlinear functions of the state, control, time, and parameters. Conversely, the \ac{ps} operators that define the transcription are linear, and their linearity can be exploited \cite{huneker2015_SPARTANImprovedGlobal}: \ac{raptor} combines exact derivatives of the linear \ac{ps} operators with numerically differentiated derivatives of the nonlinear problem-specific functions, retaining the known structure and accuracy of the \ac{ps} operators while avoiding the computational costs and potential inaccuracies of numerically differentiating the entire transcription. At a given collocation node, the problem-specific functions depend only on the local state and control values, together with the applicable times and parameters. Their state and control derivatives therefore form block-diagonal Jacobian contributions, whereas derivatives with respect to phase times and parameters form block columns, as these variables affect all nodes through the time dilation factor and the problem-specific functions. These terms are combined with the exact operator blocks, i.e., the differentiation matrix of Eq.~\eqref{eq:lagrange-defects} for Lagrange phases, or the integration matrix of Eq.~\eqref{eq:birkhoff-reconstruction} for Birkhoff phases (where the dynamics equations in Eq.~\eqref{eq:birkhoff-dynamics} remain node-local), as illustrated in Fig.~\ref{fig:hybrid-jacobian-assembly} for an example problem structure. The same separation applies to the objective gradient, for which derivatives of the running cost are combined with the quadrature weights. For more details on hybrid Jacobian computations, readers are referred to Refs.~\cite{huneker2015_SPARTANImprovedGlobal, sagliano2017_RadauPseudospectralMethod}.

\par Once the transcribed \ac{nlp} has converged, the continuous-time state is independently propagated from the optimized initial condition using the optimized controls and the known dynamics. The propagated trajectory is then compared with the interpolated state to verify dynamic consistency \cite{proulx2024_ImplementationsUniversalBirkhoff}. The solver multipliers are first returned to physical dimensions through Eq.~\eqref{eq:dual-unscaling} and then transformed into continuous-time costate and constraint-covector approximations using the covector-mapping theorem for the selected \ac{ps} method \cite{ross2001_PseudospectralTransformationCovectors, gong2008_ConnectionsCovectorMapping, ross2024_UniversalBirkhoffTheory}.

%% file: sections/results.tex
\section{Numerical Results}
\label{sec:results}

\par We now introduce four targeted benchmark problems to assess and verify the theoretical and numerical contributions presented above: a one-dimensional landing problem examines non-affine time mapping and the direct treatment of thrust-rate limits, a planar orbit-raising problem compares the conditioning and dense-grid behavior of the Lagrange and Birkhoff formulations, a dimensional atmospheric-entry problem evaluates individual and composed scaling strategies, and a five-phase \ac{rlv} problem demonstrates multiphase capabilities in a branching mission with simultaneous first-stage recovery and upper-stage insertion. Computations were performed in MATLAB R2023b on an AMD Ryzen 5 5600H processor with 16~GB of DDR4 memory. The \ac{raptor} source code is available at \url{https://github.com/raptor-solver/raptor-matlab.git}. While the \ac{raptor} framework interfaces with \ac{nlp} solvers \ac{ipopt} \cite{wachter2006_ImplementationInteriorpointFilter}, \ac{snopt} \cite{gill2002_SnoptSqpAlgorithm}, and \ac{worhp} \cite{buskens2013_ESANLPSolver}, each results subsection adopts the solver that converged more reliably in preliminary runs. In all cases, the \ac{nlp} tolerance is set to $10^{-8}$, according to the convergence criteria defined by each solver. Within each benchmark, all compared configurations use the same solver and settings.

\subsection{Non-Affine Time Mapping and Differential Constraints}
\label{subsec:landing-problem}

\par The first benchmark is a normalized one-dimensional soft-landing problem \cite{proulx2024_ImplementationsUniversalBirkhoff}, solved using \ac{ipopt}. The state is $\mathbf{x}=[h,v,m]^T$, where $h$, $v$, and $m$ denote altitude, vertical velocity, and mass, respectively, while the thrust $T$ is the control. With the gravitational acceleration normalized to unity, the fuel-optimal problem is
\begin{equation}
    \begin{aligned}
        \underset{T(\cdot),t_f}{\operatorname{minimize}}\quad
        &J=\int_{0}^{t_f}\frac{T(t)}{2.349}\,\mathrm{d}t\\
        \operatorname{subject\ to}\quad
        &\dot{h}=v\qquad
        \dot{v}=-1+\frac{T}{m}\qquad
        \dot{m}=-\frac{T}{2.349}\\
        &0\leq T\leq1.227
    \end{aligned}
    \label{eq:landing-problem}
\end{equation}
with boundary conditions
\begin{equation}
    \bigl(h(0),v(0),m(0)\bigr)=(1,-0.783,1)
    \qquad
    \bigl(h(t_f),v(t_f)\bigr)=(0,0)
    \label{eq:landing-boundary-conditions}
\end{equation}
The nominal problem produces a bang-bang thrust history \cite{proulx2024_ImplementationsUniversalBirkhoff}.

\subsubsection{Non-Affine Time Mapping}
\label{subsubsec:landing-time-mapping}

\par
\begin{figure}[htb]
    \centering
    \includegraphics[width=3.25in]{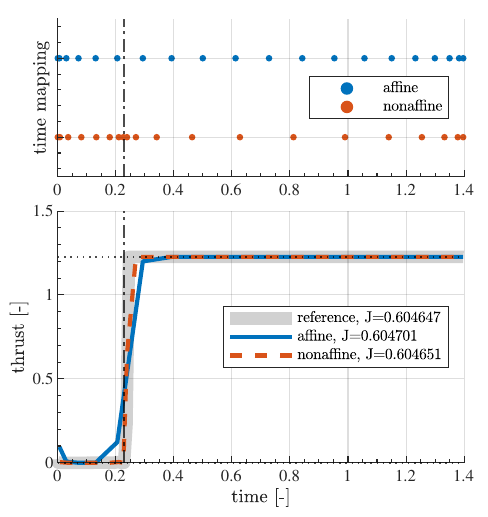}
    \caption{Affine and sinusoidal time mappings for the soft-landing problem: time-node locations (top) and thrust histories (bottom).}
    \label{fig:landing-time-dilation-effects}
\end{figure}

 To independently assess the time-mapping capability of Sec.~\ref{subsec:non-affine-time-mappings}, the nominal problem is first solved without a finite thrust-rate bound. Defining the normalized pseudotime variable $\sigma=(1+\tau)/2$, the normalized mapping function is selected as
\begin{equation}
    s(\tau;a_t,c_t)
    =
    \sigma-\frac{a_t}{2\pi}
    \left[
        \sin\!\left(2\pi(\sigma-c_t)\right)
        +
        \sin(2\pi c_t)
    \right]
    \label{eq:sinusoidal-time-map}
\end{equation}
which gives the time-dilation factor
\begin{equation}
    \rho(\tau;t_0,t_f,a_t,c_t)
    =
    \frac{t_f-t_0}{2}
    \left[
        1-a_t\cos\!\left(2\pi(\sigma-c_t)\right)
    \right]
    \label{eq:sinusoidal-time-dilation}
\end{equation}
The parameter $c_t$ locates the minimum dilation, and therefore the densest physical-time nodes, at the normalized pseudotime $\sigma=c_t$, while $a_t$ controls the concentration of the nodes around that location. The map remains monotonic for $0\leq a_t<1$ and recovers the affine map when $a_t=0$. The coarse affine solution shows a thrust switch at $t_{\mathrm{sw}}\approx0.23$, therefore the non-affine map is designed to cluster its nodes around this switch, without using the reference solution. After fixing $a_t=0.9$ (yielding a strong node concentration while remaining below the monotonicity limit $a_t<1$), $c_t$ is selected so that the corresponding physical time obtained from Eq.~\eqref{eq:normalized-time-mapping} with $t_0=0$, $t_c=t_f\,s(c_t)$, matches the estimated switch, i.e., $c_t-a_t\sin(2\pi c_t)/(2\pi)=t_{\mathrm{sw}}/t_f\approx0.164$, which yields $c_t=0.3$ for the results in Fig.~\ref{fig:landing-time-dilation-effects}. All results are generated using a Birkhoff \ac{flgr} transcription: we use the affine map with a 200-node grid to generate the high-resolution reference, and compare the performance of the affine and non-affine maps when using coarse 20-node grids.

\par Figure~\ref{fig:landing-time-dilation-effects} shows that redistributing the physical-time nodes toward the thrust-switching event improves its representation on the coarse grid: the non-affine thrust profile at the bottom of Fig.~\ref{fig:landing-time-dilation-effects} matches the reference more closely, and its fuel consumption is one order of magnitude closer to the reference than that of the affine profile. The top of Fig.~\ref{fig:landing-time-dilation-effects} compares the placement of the discretization nodes in the physical-time interval $[t_0, t_f]$, highlighting how the time mapping in Eq.~\eqref{eq:sinusoidal-time-map} aggregates time nodes around $t_c$ (vertical line).

\subsubsection{Direct Enforcement of Thrust-Rate Limits}
\label{subsubsec:landing-differential-constraints}

\par The direct treatment of differential constraints of Sec.~\ref{sec:differential-constraints} is assessed separately by imposing the thrust-rate bound $-\dot{T}_{\max}\leq\dot{T}\leq\dot{T}_{\max}$, with $\dot{T}_{\max}=4$, on the landing problem of Eq.~\eqref{eq:landing-problem}. Results are obtained using a Birkhoff \ac{flgr} transcription with 200 nodes and using affine time dilation. The problem is solved using both the direct approach of Sec.~\ref{sec:differential-constraints} and state augmentation.

\begin{figure}[htb]
    \centering
    \includegraphics[width=3.25in]{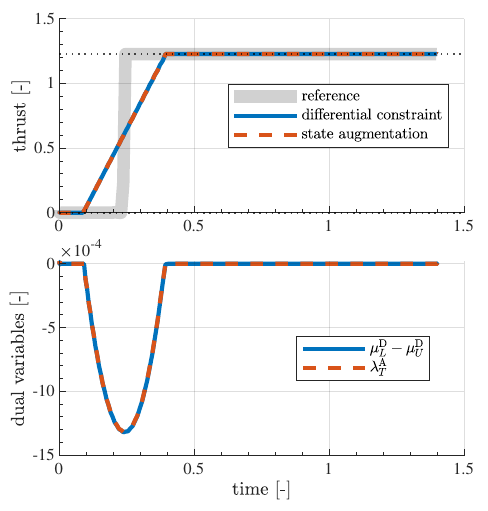}
    \caption{Direct differential constraint versus state augmentation for the rate-limited landing: thrust histories (top) and rate-constraint dual variables (bottom).}
    \label{fig:landing-primal-dual-comparison}
\end{figure}

\par The top of Fig.~\ref{fig:landing-primal-dual-comparison} shows that the two formulations produce matching thrust histories, as expected from the primal equivalence established in Sec.~\ref{subsec:control-rate-primal-equivalence}, and that the rate bound replaces the instantaneous switch of the reference solution without thrust-rate bound (Sec.~\ref{subsubsec:landing-time-mapping}) with a finite-rate thrust ramp. The bottom of Fig.~\ref{fig:landing-primal-dual-comparison} evaluates the relation in Eq.~\eqref{eq:implicit-control-costate} by comparing the multiplier difference $\mu_L^{\mathrm{D}}-\mu_U^{\mathrm{D}}$ of the direct formulation with the costate $\lambda_T^{\mathrm{A}}$ of the augmented thrust state, where the superscripts D and A denote the direct and state-augmentation formulations, respectively. The multiplier difference departs from zero on the thrust-increasing arc, where the upper rate bound is active, and vanishes where neither rate bound is active. Its agreement with the augmented-state costate $\lambda_T^{\mathrm{A}}$ numerically supports the dual equivalence established in Sec.~\ref{subsec:control-rate-dual-equivalence}.

\begin{table}[htb]
    \centering
    \caption{Direct differential constraint versus state augmentation for the rate-limited landing: transcription structure and solver performance.}
    \label{tab:landing-transcription-comparison}
    \setlength{\tabcolsep}{5pt}
    \begin{tabular}{@{}lcc@{}}
        \hline\hline
        Metric & Direct & Augmented \\
        \hline
        Jacobian size, $n_c\times n_z$ & $1598 \times 1602$ & $1598 \times 1802$ \\
        Jacobian nonzeros & 163186 & 163385 \\
        \ac{ipopt} iterations & 27 & 35 \\
        Average solve time$^{\mathrm{a}}$ (s) & 5.649 & 8.173 \\
        \hline\hline
        \multicolumn{3}{@{}l@{}}{\footnotesize $^{\mathrm{a}}$Averaged over 100 runs on the same machine.}
    \end{tabular}
\end{table}

\par As summarized in Table~\ref{tab:landing-transcription-comparison}, state augmentation adds the 200 decision variables of the thrust-rate control and 199 Jacobian nonzeros, consistent with the additional collocation constraints for the augmented thrust state. The constraint count is unchanged since the augmented rate limits are imposed as variable bounds instead of the nonlinear inequalities of the direct formulation. For this benchmark, state augmentation also requires more \ac{ipopt} iterations and a longer average solve time.

\subsection{Dense-Grid Lagrange-Birkhoff Comparison}
\label{subsec:orbit-raising-problem}

\par The second benchmark examines the dense-grid behavior of the Lagrange and Birkhoff formulations using a normalized planar low-thrust orbit-raising problem based on Ref.~\cite{longuski2014_OptimalControlAerospace}, solved using \ac{snopt}. The state is $\mathbf{x}=[r,\theta,v_r,v_\theta]^T$ and includes orbital radius, longitude, and radial and tangential velocity, while the thrust angle $\alpha$ is the control. The gravitational parameter is normalized to unity, and the constant thrust acceleration is $T=4\times10^{-4}$. For a prescribed final time, the terminal radius is maximized according to
\begin{equation}
    \begin{aligned}
        \underset{\alpha(\cdot)}{\operatorname{maximize}}\quad
        &J=r(t_f)\\
        \operatorname{subject\ to}\quad
        &\dot{r}=v_r\qquad
        \dot{\theta}=\frac{v_\theta}{r}\\
        &\dot{v}_r=\frac{v_\theta^2}{r}-\frac{1}{r^2}+T\sin\alpha \qquad
        \dot{v}_\theta=-\frac{v_rv_\theta}{r}+T\cos\alpha
    \end{aligned}
    \label{eq:orbit-raising-problem}
\end{equation}
with boundary conditions
\begin{equation}
    t_0=0\qquad t_f=40\pi\qquad
    \mathbf{x}(0)=
    \begin{bmatrix}
        1 & 0 & 0 & 1
    \end{bmatrix}^{T}
    \qquad
    v_r(t_f)=0
    \qquad
    v_\theta(t_f)=\frac{1}{\sqrt{r(t_f)}}
    \label{eq:orbit-raising-boundary-conditions}
\end{equation}
The terminal longitude is free, the fixed duration corresponds to 20 periods of the initial circular orbit, and $-\pi\leq\alpha\leq\pi$. The multiple revolutions and oscillatory steering history make this problem a representative benchmark for dense-grid methods: we transcribed it using both the Birkhoff and Lagrange formulations with node counts varied over $N\in\{25,50,100,150,\ldots,500\}$, while the rest of the problem definition is held constant. The independently propagated trajectory is obtained using a variable-order (4-5) Runge-Kutta integrator, following the validation procedure described in Sec.~\ref{subsec:hybrid-jacobians}.

\par
\begin{figure*}[htb]
    \centering
    \begin{subfigure}{3.25in}
        \centering
        \includegraphics[width=\linewidth]{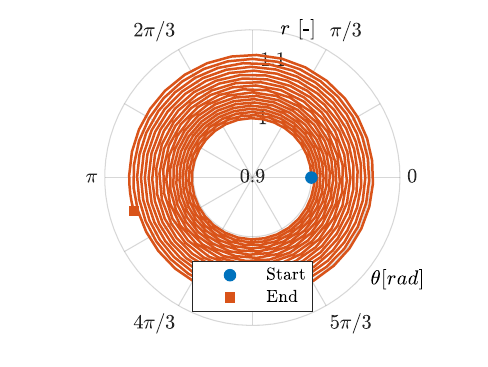}
        \caption{}
        \label{subfig:orbit-raising-trajectory}
    \end{subfigure}%
    \hfill
    \begin{subfigure}{3.25in}
        \centering
        \includegraphics[width=\linewidth]{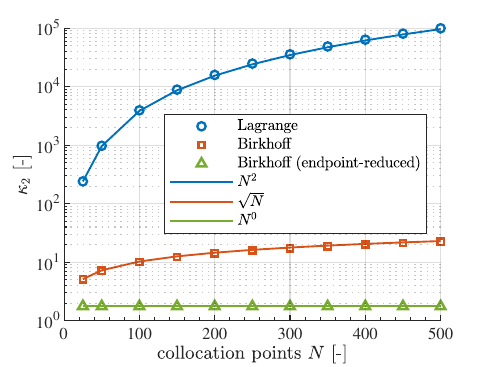}
        \caption{}
        \label{subfig:orbit-raising-conditioning}
    \end{subfigure}
    \par\medskip
    \begin{subfigure}{3.25in}
        \centering
        \includegraphics[width=\linewidth]{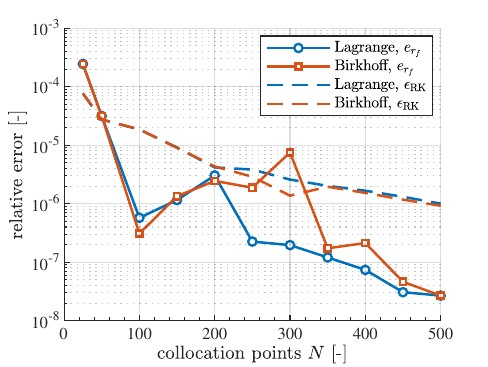}
        \caption{}
        \label{subfig:orbit-raising-accuracy}
    \end{subfigure}%
    \hfill
    \begin{subfigure}{3.25in}
        \centering
        \includegraphics[width=\linewidth]{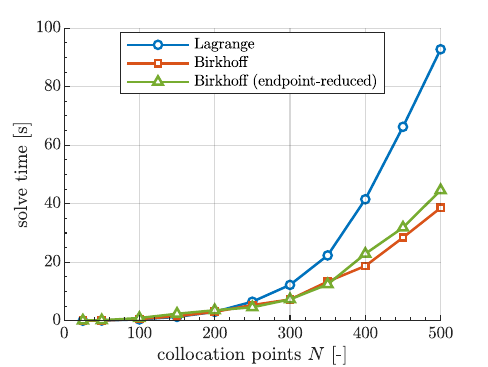}
        \caption{}
        \label{subfig:orbit-raising-solve-time}
    \end{subfigure}
    \caption{Lagrange-Birkhoff comparison for orbit raising: a) representative trajectory; b) operator conditioning; c) accuracy and propagation error; d) \ac{snopt} solve time.}
    \label{fig:orbit-raising-comparison}
\end{figure*}

 Figure~\ref{fig:orbit-raising-comparison} summarizes the performance of the tested formulations across the increasing node counts. The conditioning comparison uses the trajectory-independent linear collocation systems rather than the complete nonlinear dynamics Jacobian, thereby isolating the numerical properties of the interpolation operators from the trajectory, dynamics, and scaling. For the Lagrange formulation, the metric is $\kappa_2(\mathbf{D})$. For the Birkhoff formulation, the reported metric is the $\kappa_2(\mathbf{C}^a)$ used in Sec.~\ref{subsec:new-birkhoff-interpolants}. Accuracy is evaluated through two metrics: the relative terminal-radius error $e_{r_f}$, which compares the maximized final radius against the finest-mesh ($N=500$) solution of the same formulation, and the relative root-mean-square state error $\epsilon_{\mathrm{RK}}$, which compares the state history with the independently propagated Runge-Kutta trajectory. Solve times are averaged over 100 identical runs for each configuration. Conditioning and solve times are distinguished between the regular and endpoint-reduced Birkhoff \ac{flgr} approaches, while their accuracy is not distinguished, as it would only report redundant data.

\par Figure~\ref{subfig:orbit-raising-trajectory} provides the geometric context for the dense-grid comparison. Figure~\ref{subfig:orbit-raising-conditioning} clearly shows the expected trends in \ac{ps}-specific condition numbers: the Lagrange formulation closely follows $\mathcal{O}(N^2)$, the Birkhoff one shows the expected $\mathcal{O}(\sqrt{N})$ behavior, and the endpoint-reduced version exhibits the bounded flat $\mathcal{O}(1)$ trend. Figure~\ref{subfig:orbit-raising-accuracy} expectedly shows decreasing error trends for refined grids. As expected from the connection between the Lagrange and Birkhoff formulations drawn in Sec.~\ref{subsubsec:flgr-lagrange-connection}, the two formulations exhibit similar accuracy, while their differences mainly concern solver performance. Figure~\ref{subfig:orbit-raising-solve-time} reports the corresponding \ac{snopt} solve times. We notice that Birkhoff formulations are generally more computationally convenient for increasingly dense grids, despite the lifted representation via the virtual derivative variables. No significant difference can be identified between the \ac{flgr} Birkhoff formulation and its endpoint reduced counterpart: a flat condition number of the \ac{ps}-specific portion of the transcription does not necessarily translate into improved performance.

\subsection{Composable Scaling for Atmospheric Entry}
\label{subsec:entry-problem}

\par To isolate the effect of variable and constraint scaling, the dimensional Space Shuttle atmospheric-entry benchmark from Ref.~\cite{betts2010_PracticalMethodsOptimal} is considered, solved using \ac{worhp}. This problem defines states $\mathbf{x}=[h,\phi,\theta,v,\gamma,\psi]^T$, angle of attack and bank angle controls $\mathbf{u}=[\alpha,\beta]^T$, and the final latitude $\theta$ is maximized subject to
\begin{equation}
    \begin{aligned}
        \dot h&=v\sin\gamma
        &\dot\theta&=\frac{v}{r}\cos\gamma\cos\psi\\
        \dot\phi&=\frac{v}{r}\frac{\cos\gamma\sin\psi}{\cos\theta}
        &\dot v&=-\frac{D}{m}-g\sin\gamma\\
        \dot\gamma&=\frac{L}{mv}\cos\beta+\cos\gamma\left(\frac{v}{r}-\frac{g}{v}\right)&
        \dot\psi&=\frac{L\sin\beta}{mv\cos\gamma}
        +\frac{v}{r}\cos\gamma\sin\psi\tan\theta\\
        r&=R_e+h\quad g=\frac{\mu}{r^2}\quad
        \rho_a=\rho_0e^{-h/h_r}&
        (D,L)&=\frac{1}{2}\rho_a Sv^2(c_D,c_L)\\
        q&=q_a\,17700\sqrt{\rho_a}(10^{-4}v)^{3.07}\leq q_U
    \end{aligned}
    \label{eq:entry-problem}
\end{equation}
where $\hat{\alpha}=180\alpha/\pi$ is the angle of attack in degrees and the model coefficients are computed as
\begin{equation}
    c_L=a_0+a_1\hat{\alpha},\qquad
    c_D=b_0+b_1\hat{\alpha}+b_2\hat{\alpha}^2,\qquad
    q_a=c_0+c_1\hat{\alpha}+c_2\hat{\alpha}^2+c_3\hat{\alpha}^3.
    \label{eq:entry-coefficients}
\end{equation}
The model parameters are retained from Ref.~\cite{betts2010_PracticalMethodsOptimal} and summarized using the US customary units from the original source in Table~\ref{tab:entry-parameters}.
\begin{table*}[htb]
    \centering
    \caption{Entry-problem vehicle, atmospheric, and aerodynamic parameters (US customary units).}
    \label{tab:entry-parameters}
    \setlength{\tabcolsep}{8pt}
    \begin{tabular}{@{}cccccc@{}}
        \hline\hline
        Parameter & Value & Parameter & Value & Parameter & Value \\
        \hline
        $R_e$     & $2.09\times10^{7}$ 
        & $\mu$   & $1.41\times10^{16}$ 
        & $\rho_0$ & $2.38\times10^{-3}$ \\

        $h_r$     & $2.38\times10^{4}$ 
        & $S$     & $2690$ 
        & $m$     & $6.309\times10^{3}$ \\

        $g_0$     & $32.174$ 
        & $q_U$        & $70$
        &         & \\
        \hline
        $a_0$     & $-0.20704$ 
        & $a_1$   & $0.029244$ 
        & $b_0$   & $0.07854$ \\

        $b_1$     & $-6.16\times10^{-3}$ 
        & $b_2$   & $6.21\times10^{-4}$ 
        & $c_0$   & $1.0672$ \\

        $c_1$     & $-1.92\times10^{-2}$ 
        & $c_2$   & $2.13\times10^{-4}$ 
        & $c_3$   & $-1.01\times10^{-6}$ \\
        \hline\hline
    \end{tabular}
\end{table*}
The free-final-time boundary data and relevant box bounds are compactly summarized as
\begin{equation}
    \begin{gathered}
        \mathbf{x}(0)=
        [260000,0,0,25600,-1^\circ,90^\circ]^T
        \qquad
        [h(t_f),v(t_f),\gamma(t_f)]^T=[80000,2500,-5^\circ]^T\\
        h\geq0\quad v\geq1\quad
        |\theta|,|\gamma|\leq89^\circ\quad
        |\alpha|\leq90^\circ\quad
        -89^\circ\leq\beta\leq1^\circ
    \end{gathered}
    \label{eq:entry-boundary-conditions}
\end{equation}
Here altitude and velocity are expressed in feet and feet per second, respectively, while all angle variables are expressed in radians. The dimensional spread and highly nonlinear dynamics of this problem make it a relevant test case for the \ac{raptor} composable scaling pipeline.

\par
\begin{figure*}[htb]
    \centering
    \begin{subfigure}{3.25in}
        \centering
        \includegraphics[width=\linewidth]{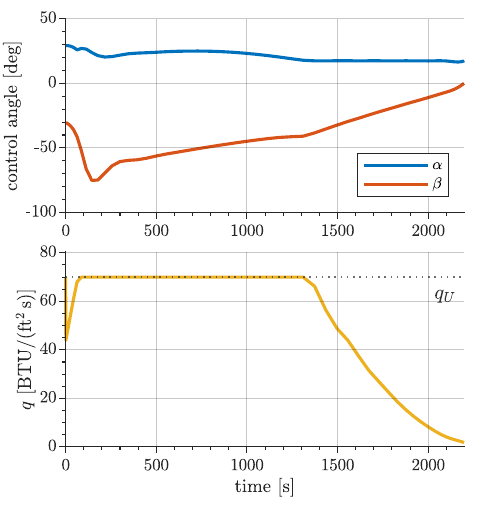}
        \caption{}
        \label{subfig:entry-controls-heating}
    \end{subfigure}%
    \hfill
    \begin{subfigure}{3.25in}
        \centering
        \includegraphics[width=\linewidth]{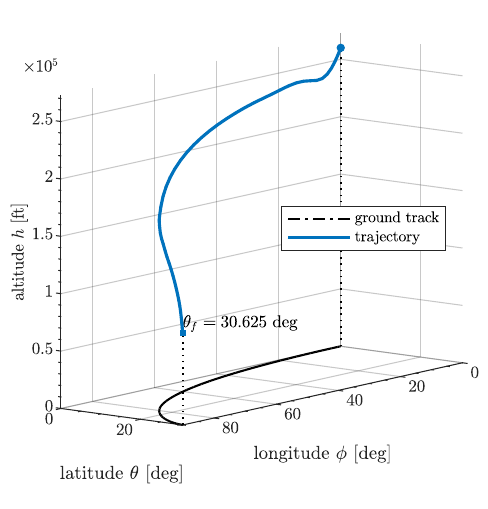}
        \caption{}
        \label{subfig:entry-trajectory}
    \end{subfigure}
    \caption{Shuttle entry solution: a) control angles (top) and heat flux (bottom); b) trajectory and ground track.}
    \label{fig:entry-solution-summary}
\end{figure*}

 The problem is collocated on a Birkhoff \ac{flgr} grid with 100 nodes. An example \ac{raptor} solution (obtained using \ac{pjrn} scaling) is concisely summarized in Fig.~\ref{fig:entry-solution-summary}, using the US customary units of Ref.~\cite{betts2010_PracticalMethodsOptimal} to ease the comparison. The control angles at the top of Fig.~\ref{subfig:entry-controls-heating} closely match those from Ref.~\cite{betts2010_PracticalMethodsOptimal}, and the heat flux at the bottom quickly reaches and maintains its upper bound value $q_U$. Figure~\ref{subfig:entry-trajectory} shows the corresponding trajectory in the longitude-latitude-altitude space, together with its ground track, whose final latitude $\theta(t_f)$ is the maximized crossrange objective: both $\theta(t_f)=30.625$~deg and $t_f=2198.7$~s match the results obtained in Ref.~\cite{betts2010_PracticalMethodsOptimal}.

\par The complete constraint Jacobian is evaluated at the same physical initial guess using the scaling approaches introduced in Sec.~\ref{subsec:composable-scaling} and their ordered compositions. The comparison therefore uses the effective condition number $\kappa_2^+=\sigma_1/\sigma_r$, where $\sigma_r$ is the smallest singular value satisfying $\sigma_r>\epsilon_{\mathrm{rank}}$. For consistent comparisons across all configurations, we define the numerical-rank tolerance as $\epsilon_{\mathrm{rank}}=\max(n_c,n_z)\operatorname{eps}(\sigma_1)$, where $n_c$ and $n_z$ are the numbers of constraints and decision variables, respectively, and $\operatorname{eps}(\sigma_1)$ is the current machine precision multiplied by $|\sigma_1|$.

\begin{table*}[htb]
    \centering
    \caption{Effective initial-guess Jacobian conditioning and mean solver performance for the atmospheric-entry problem; bold values mark the best results.}
    \label{tab:entry-scaling-comparison}
    \setlength{\tabcolsep}{8pt}
    \begin{tabular}{@{}lccc@{}}
        \hline\hline
        Scaling sequence & $\kappa_2^+(\mathbf{J}_{c,0})$ & Mean solve time (s) & $\theta(t_f)$ (deg) \\
        \hline
        0: None                              & $1.09\times10^{13}$ & --- & --- \\
        1: Bounds projection                 & $1.20\times10^{9}$  & ---  & --- \\
        2: \ac{jrn}                          & $4.91\times10^{7}$  & 55.126  & \textbf{30.625} \\
        3: \ac{pjrn}                         & $3.71\times10^{3}$  & 4.516   & \textbf{30.625} \\
        4: Static Ruiz                       & $1.10\times10^{4}$  & ---  & --- \\
        5: Bounds projection $+$ static Ruiz & $7.58\times10^{3}$  & ---  & --- \\
        6: \ac{jrn} $+$ static Ruiz          & {\boldmath$2.15\times10^{3}$}  & 10.124  & 30.217 \\
        7: Static Ruiz $+$ \ac{jrn}          & $7.90\times10^{3}$  & 14.791  & 30.233 \\
        8: Bounds projection $+$ static Ruiz $+$ \ac{jrn}
                                          & $6.73\times10^{3}$  & 6.387   & \textbf{30.625} \\
        9: \ac{pjrn} $+$ static Ruiz         & $2.22\times10^{3}$  & \textbf{2.283}   & \textbf{30.625} \\
        \hline\hline
    \end{tabular}
\end{table*}

\par Table~\ref{tab:entry-scaling-comparison} reports the effective initial-guess conditioning and the \ac{worhp} performance for each scaling configuration, averaged over 100 runs, and dashes denote non-converged configurations. All converged configurations include \ac{jrn} (cases 2, 3, and 6--9), whereas cases 0, 1, 4, and 5 fail, even though static Ruiz equilibration alone reduces $\kappa_2^+$ by approximately nine orders of magnitude. Combining bounds projection with \ac{jrn} (cases 3, 8, and 9) consistently recovers the reference final latitude of 30.625~deg, and reduces the mean solve time from 55.1~s (case 2) to 4.5~s (case 3). Without bounds projection, combining \ac{jrn} with static Ruiz equilibration (cases 6 and 7) yields a lower final latitude of approximately 30.2~deg, despite case 6 exhibiting the lowest $\kappa_2^+$ among all configurations. The order of the transformations also matters, as anticipated by Eqs.~\eqref{eq:composed-variable-scaling} and~\eqref{eq:composed-function-scaling}: cases 6 and 7, as well as cases 8 and 9, combine the same transformations in different orders and differ in both conditioning and solve time. In particular, applying static Ruiz equilibration after \ac{pjrn} (case 9) reduces $\kappa_2^+$ from $3.71\times10^3$ to $2.22\times10^3$ and halves the mean solve time, while recovering the same final latitude. Overall, the effective condition number alone does not predict solver performance, and no single strategy is uniformly preferable: the composition mechanism of Sec.~\ref{subsec:composable-scaling} allows problem-tailored sequences to be quickly assembled and compared without modifying the transcription, enabling fast iterations to identify effective scaling methods.

\subsection{Heterogeneous Multiphase Launch-Vehicle Demonstration}
\label{subsec:rlv-problem}

\par The final benchmark is a branching launch-vehicle mission, consisting of a shared first-stage ascent coupled at stage separation to simultaneous booster \ac{rtls} and upper-stage insertion trajectories, and it is used to demonstrate the developed features within a single multiphase problem. The mission is represented using a planar, nonrotating, flat-Earth with constant gravitational acceleration and an exponential atmosphere model. For each phase $q$, the state and control vectors are
\begin{equation}
    \mathbf{x}^{(q)}=
    \begin{bmatrix}
        s & h & v_s & v_h & m
    \end{bmatrix}^{T}
    \qquad
    \mathbf{u}^{(q)}=
    \begin{bmatrix}
        \eta & \theta
    \end{bmatrix}^{T}
    \label{eq:rlv-state-control}
\end{equation}
where $s$ is the downrange position, $h$ is the altitude, $v_s$ and $v_h$ are the corresponding velocity components, $m$ is the vehicle mass, $\eta$ is the throttle fraction, and $\theta$ is the thrust-direction angle measured counterclockwise from the positive downrange axis.

\par The atmospheric density, vehicle speed, and drag-force components are defined as
\begin{equation}
    \rho(h)=\rho_0\exp\left(-\frac{h}{H}\right)
    \qquad
    v=\sqrt{v_s^2+v_h^2}
    \qquad
    \begin{bmatrix}
        F_{D,s}\\
        F_{D,h}
    \end{bmatrix}
    =
    -\frac{1}{2}\rho(h)C_DAv
    \begin{bmatrix}
        v_s\\
        v_h
    \end{bmatrix}
    \label{eq:rlv-atmosphere-drag}
\end{equation}
For a powered phase, we define $T_{\max}^{(q)}$ and $I_{sp}^{(q)}$ as the constant phase-effective thrust and specific impulse, respectively. The equations of motion are
\begin{align}
    \dot{s}
    &=
    v_s
    &
    \dot{h}
    &=
    v_h
    \label{eq:rlv-kinematics}\\
    \dot{v}_s
    &=
    \frac{\eta T_{\max}^{(q)}\cos\theta+F_{D,s}}{m}
    &
    \dot{v}_h
    &=
    \frac{\eta T_{\max}^{(q)}\sin\theta+F_{D,h}}{m}-g_0
    \label{eq:rlv-translational-dynamics}\\
    \dot{m}
    &=
    -\frac{\eta T_{\max}^{(q)}}{g_0I_{sp}^{(q)}}
    &&
    \label{eq:rlv-mass-dynamics}
\end{align}
During the unpowered free-fall phase, $\eta=0$ and the vehicle mass remains constant. The propulsion model uses constant effective thrust levels and omits altitude-dependent performance, throttle envelopes, and engine-to-engine variations. The first-stage values preserve the assumed engine-count scaling of $T_{\max}^{(\mathrm{A})}:T_{\max}^{(\mathrm{B})}:T_{\max}^{(\mathrm{L})}=9:3:1$ using a nominal effective thrust of $654\ \mathrm{kN}$ per active engine (see Table~\ref{tab:rlv-parameters}). This model captures the propulsion changes associated with staging and recovery without additional engine-performance states or maps.

\par The principal vehicle and environmental parameters are summarized in Table~\ref{tab:rlv-parameters}. The target horizontal velocity $v_{\mathrm{orb}}$ is computed as the circular orbit velocity corresponding to the target altitude $h_{\mathrm{orb}}$ through $v_{\mathrm{orb}}=\sqrt{\frac{\mu_E}{R_E+h_{\mathrm{orb}}}}=7.67259\ \mathrm{km/s}$. Because Earth curvature and inverse-square gravity are not included in Eqs.~\eqref{eq:rlv-kinematics}--\eqref{eq:rlv-mass-dynamics}, the $v_{\mathrm{orb}}$ target is used as a terminal insertion surrogate rather than as part of a high-fidelity orbital model.

\begin{table*}[htb]
    \centering
    \caption{Environmental and Falcon~9-inspired vehicle parameters.}
    \label{tab:rlv-parameters}
    \setlength{\tabcolsep}{4pt}
    \begin{tabular}{@{}lll@{\hspace{1.2em}}lll@{}}
        \hline\hline
        Parameter & Symbol & Value & Parameter & Symbol & Value \\
        \hline
        Gravity & $g_0$ & $9.80665\ \mathrm{m/s^2}$ & Earth radius & $R_E$ & $6371\ \mathrm{km}$ \\
        Sea-level density & $\rho_0$ & $1.225\ \mathrm{kg/m^3}$ & Earth grav. parameter & $\mu_E$ & $3.986\times10^{14}\ \mathrm{m^3/s^2}$ \\
        Scale height & $H$ & $7.5\ \mathrm{km}$ & Drag coefficient & $C_D$ & $0.75$ \\
        Reference area & $A$ & $10.5209\ \mathrm{m^2}$ & Target altitude & $h_{\mathrm{orb}}$ & $400\ \mathrm{km}$ \\
        Stage-1 dry mass & $m_{1,d}$ & $25600\ \mathrm{kg}$ & Stage-1 propellant & $m_{1,p}$ & $395700\ \mathrm{kg}$ \\
        Stage-2 dry mass & $m_{2,d}$ & $4500\ \mathrm{kg}$ & Stage-2 propellant & $m_{2,p}$ & $107000\ \mathrm{kg}$ \\
        Payload mass & $m_{\mathrm{PL}}$ & $15000\ \mathrm{kg}$ & Stage-1 specific impulse & $I_{sp,1}$ & $283\ \mathrm{s}$ \\
        Stage-2 specific impulse & $I_{sp,2}$ & $348\ \mathrm{s}$ & Ascent thrust & $T_{\max}^{(\mathrm{A})}$ & $5.886\ \mathrm{MN}$ \\
        Boostback thrust & $T_{\max}^{(\mathrm{B})}$ & $1.962\ \mathrm{MN}$ & Landing thrust & $T_{\max}^{(\mathrm{L})}$ & $0.654\ \mathrm{MN}$ \\
        Stage-2 thrust & $T_{\max}^{(\mathrm{S})}$ & $0.981\ \mathrm{MN}$ & Target speed & $v_{\mathrm{orb}}$ & $7.67259\ \mathrm{km/s}$ \\
        Initial mass & $m_0$ & $547800\ \mathrm{kg}$ & Stage-2 initial mass & $m_{2, 0}$ & $126500\ \mathrm{kg}$ \\
        \hline\hline
    \end{tabular}
\end{table*}

\par The trajectory comprises ascent $\mathrm{A}$, boostback $\mathrm{B}$, free fall $\mathrm{F}$, landing $\mathrm{L}$, and second-stage ascent $\mathrm{S}$. The combined vehicle ascends to stage separation, where the trajectory branches into the first-stage return sequence $\mathrm{B}\rightarrow\mathrm{F}\rightarrow\mathrm{L}$ and the upper-stage phase $\mathrm{S}$. All phases are collocated on \ac{flgr} grids, while the interpolant is selected phase by phase exploiting the heterogeneous transcription of Sec.~\ref{subsec:heterogeneous-multiphase-transcription}. The Birkhoff formulation is adopted for the long and densely discretized ascent and second-stage arcs, where the conditioning advantage observed in Sec.~\ref{subsec:orbit-raising-problem} is most relevant. The Lagrange formulation is instead adopted for the short, low-order boostback, free-fall, and landing phases, where the conditioning of the differentiation matrix remains moderate and the virtual derivative variables of the Birkhoff formulation would only enlarge the \ac{nlp}. Table~\ref{tab:rlv-phases} summarizes the transcription, collocation-point counts, control restrictions, terminal-time guesses and related bounds for each phase.

\begin{table*}[htb]
    \centering
    \caption{Phase definitions and discretization of the \ac{rlv} benchmark; the last column reports the terminal-time initial guess and its admissible interval.}
    \label{tab:rlv-phases}
    \setlength{\tabcolsep}{4pt}
    \begin{tabular}{@{}llllllll@{}}
        \hline\hline
        Phase & Mode & Grid & Interpolant & $N_{c}$ & $\eta$ & $\theta$ & $t_f^0\;[t_{f,L},t_{f,U}]$ (s) \\
        \hline
        Ascent, $\mathrm{A}$ & Combined, 9 engines & \ac{flgr} & Birkhoff & $100$ & $1$ & $[0,\pi]$ & $140\;[90,190]$ \\
        Boostback, $\mathrm{B}$ & Stage 1, 3 engines & \ac{flgr} & Lagrange & $10$ & $1$ & $\pi$ & $181\;[131,231]$ \\
        Free fall, $\mathrm{F}$ & Stage 1, unpowered & \ac{flgr} & Lagrange & $30$ & $0$ & $0$ & $361\;[311,411]$ \\
        Landing, $\mathrm{L}$ & Stage 1, 1 engine & \ac{flgr} & Lagrange & $50$ & $1$ & $[0,\pi]$ & $451\;[401,501]$ \\
        Stage 2, $\mathrm{S}$ & Stage 2, 1 engine & \ac{flgr} & Birkhoff & $100$ & $[0,1]$ & $[0,\pi]$ & $1000\;[0,2000]$ \\
        \hline\hline
    \end{tabular}
\end{table*}

\par The combined vehicle departs from rest at the origin at $t_{0,\mathrm{A}}=0$ with launch mass $m_0$. At stage separation, both outgoing branches inherit the final ascent position and velocity, while their initial times and masses are linked as
\begin{equation}
    t_{0,\mathrm{B}}=t_{0,\mathrm{S}}=t_{f,\mathrm{A}}
    \qquad
    m^{(\mathrm{B})}(t_{0,\mathrm{B}})=m^{(\mathrm{A})}(t_{f,\mathrm{A}})-m_{2,0}
    \qquad
    m^{(\mathrm{S})}(t_{0,\mathrm{S}})=m_{2,0}
    \label{eq:rlv-separation-linkage}
\end{equation}
Full state and time continuity are imposed between the consecutive first-stage return phases $\mathrm{B}\rightarrow\mathrm{F}\rightarrow\mathrm{L}$, and the first stage is required to return to the launch site ($s=h=0$) with zero terminal velocity and a final mass not lower than its dry mass $m_{1,d}$. During the landing phase, the vehicle is constrained to satisfy the glide-slope constraint
\begin{equation}
    h^{(\mathrm{L})}(t)-s^{(\mathrm{L})}(t)\tan\gamma_{\mathrm{gs}}\geq0
    \qquad
    \gamma_{\mathrm{gs}}=70\ \mathrm{deg}
    \label{eq:rlv-landing-glideslope}
\end{equation}

\par The upper-stage terminal altitude and velocity are constrained to the surrogate circular-orbit conditions
\begin{equation}
    h^{(\mathrm{S})}(t_{f,\mathrm{S}})=h_{\mathrm{orb}}
    \qquad
    v_s^{(\mathrm{S})}(t_{f,\mathrm{S}})=v_{\mathrm{orb}}
    \qquad
    v_h^{(\mathrm{S})}(t_{f,\mathrm{S}})=0
    \label{eq:rlv-second-stage-terminal}
\end{equation}
while its terminal downrange position and mass remain free. The objective is
\begin{equation}
    \underset{\mathbf{z}}{\operatorname{minimize}}
    \quad
    J=-m^{(\mathrm{S})}(t_{f,\mathrm{S}})
    \label{eq:rlv-objective}
\end{equation}
which maximizes the upper-stage mass delivered into the target orbit. All computations are performed in SI units, while selected results are displayed in kilometers and kilometers per second for readability. The \ac{ipopt} solver is selected for this problem.

\par
\begin{figure*}[htb]
    \centering
    \begin{subfigure}{3.25in}
        \centering
        \includegraphics[width=\linewidth]{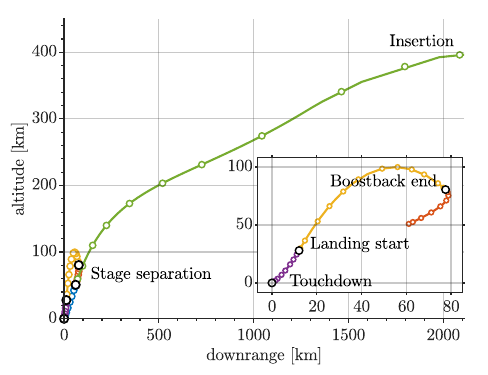}
        \caption{}
        \label{subfig:rlv-mission-geometry}
    \end{subfigure}%
    \hfill
    \begin{subfigure}{3.25in}
        \centering
        \includegraphics[width=\linewidth]{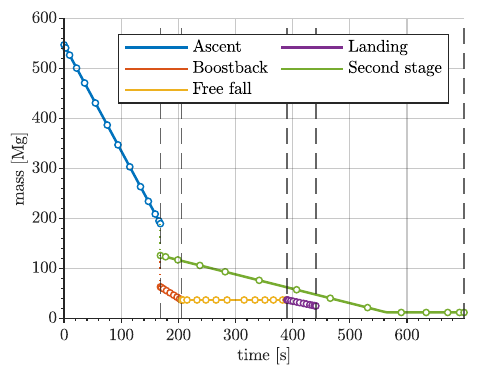}
        \caption{}
        \label{subfig:rlv-mass}
    \end{subfigure}
    \par\medskip
    \begin{subfigure}{3.25in}
        \centering
        \includegraphics[width=\linewidth]{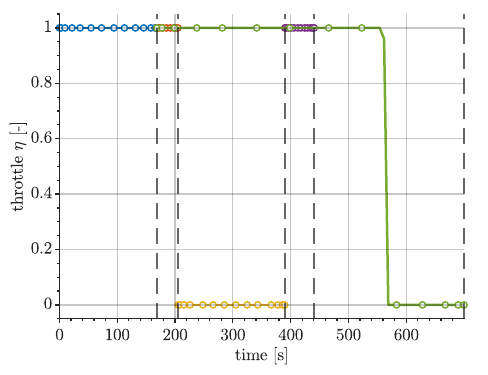}
        \caption{}
        \label{subfig:rlv-throttle}
    \end{subfigure}%
    \hfill
    \begin{subfigure}{3.25in}
        \centering
        \includegraphics[width=\linewidth]{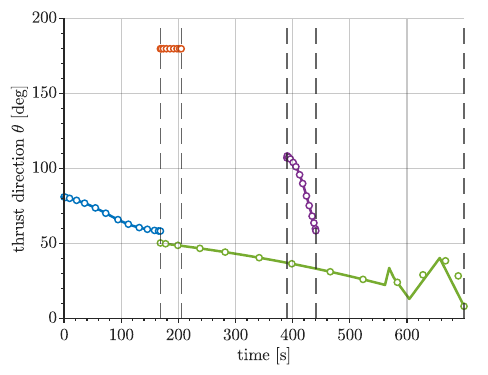}
        \caption{}
        \label{subfig:rlv-thrust-direction}
    \end{subfigure}
    \caption{Five-phase \ac{rlv} solution: a) mission geometry; b) mass histories; c) throttle histories; d) thrust-direction histories.}
    \label{fig:rlv-reference-solution}
\end{figure*}

 The resulting representative solution is independently validated via numerical propagation and summarized in Fig.~\ref{fig:rlv-reference-solution}, which reports the geometry of both mission branches (Fig.~\ref{subfig:rlv-mission-geometry}) together with their mass and control histories (Figs.~\ref{subfig:rlv-mass}--\ref{subfig:rlv-thrust-direction}). Continuous lines represent the independently propagated quantities, while markers show the nodal values at a subset of at most 14 nodes per phase, sampled uniformly in the node index to preserve the clustering of the \ac{ps} nodes.

\par The shared ascent produces a single optimized separation state from which the first stage reverses its downrange motion and returns to the launch site while the upper stage satisfies the insertion conditions. Since the first-stage throttle is prescribed phase by phase (Table~\ref{tab:rlv-phases}), the return sequence follows an imposed full-thrust/coast/full-thrust structure whose switching times are optimized through the free phase durations. The upper-stage throttle, left free within $[0,1]$, converges to a bang-bang profile consistent with maximizing its terminal mass (Fig.~\ref{subfig:rlv-throttle}), while its thrust direction remains nearly tangential to satisfy the insertion constraints with limited propellant expenditure (Fig.~\ref{subfig:rlv-thrust-direction}). After the throttle switches to zero, the thrust direction no longer enters the dynamics or the cost, thus $\theta$ is not determined by optimality conditions and the oscillations at the end of the second-stage arc in Fig.~\ref{subfig:rlv-thrust-direction} reflect its indetermination, rather than a physical feature of the solution. The agreement between the discrete solution and the independent continuous propagation supports the dynamic consistency of the coupled multiphase transcription.

\par
\begin{figure}[htb]
    \centering
    \includegraphics[width=3.25in]{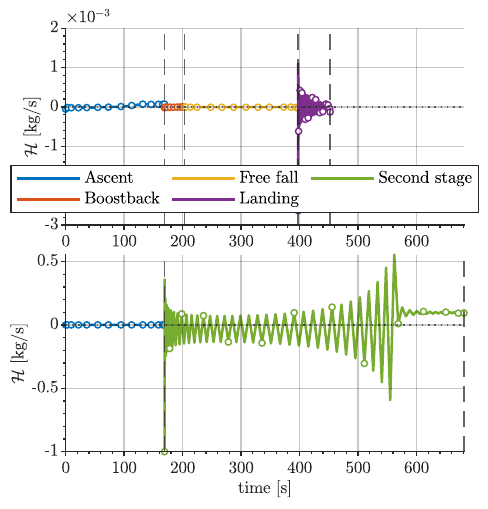}
    \caption{Phase Hamiltonians of the \ac{rlv} solution: first-stage branch (top) and second-stage branch (bottom); dashed lines mark phase interfaces.}
    \label{fig:rlv-hamiltonian}
\end{figure}

 The dual consistency of the heterogeneous transcription is assessed through the Hamiltonian of each phase, computed from the costates recovered via the covector mapping of the corresponding Lagrange \cite{ross2001_PseudospectralTransformationCovectors, gong2008_ConnectionsCovectorMapping} or Birkhoff \cite{ross2024_UniversalBirkhoffTheory} formulation. Since none of the phase dynamics depends explicitly on time, the Hamiltonian is constant along each optimal phase \cite{ross2015_PrimerPontryaginsPrinciple, longuski2014_OptimalControlAerospace}. Its behavior therefore provides a transcription-independent check of the recovered dual information, which is also sensitive to inconsistencies in the time dilation factors and phase linkages. Figure~\ref{fig:rlv-hamiltonian} reports the resulting Hamiltonian histories for the first-stage and second-stage branches.

\par Figure~\ref{fig:rlv-hamiltonian} shows that the Hamiltonian remains close to zero in all phases: since all phase durations are free and none of the phase-time bounds is active, it holds that $\mathcal{H}\equiv0$. Its magnitude is assessed with respect to the dominant terms of the Hamiltonian, e.g., the mass-flow contribution $|\lambda_m\dot{m}|\approx 287$~kg/s in the second stage. In the first-stage branch (top of Fig.~\ref{fig:rlv-hamiltonian}), $|\mathcal{H}|$ remains below $10^{-4}$~kg/s along the Birkhoff ascent and the Lagrange boostback and free-fall phases, without discontinuities at their interfaces, while localized oscillations of up to $3\times10^{-3}$~kg/s appear at the beginning of the landing phase. In the second stage (bottom of Fig.~\ref{fig:rlv-hamiltonian}), apart from an isolated value of approximately $-1$~kg/s (about 0.35\% of the dominant term) at its first node, node-to-node oscillations of up to approximately $0.6$~kg/s (about 0.2\% of the dominant term) develop along the powered arc and vanish after the throttle switch at $t\approx 565$~s, where the Hamiltonian settles at a small constant offset of approximately $0.1$~kg/s (about 0.03\% of the dominant term) instead of zero. Both effects are consistent with the approximation of the discontinuous bang-bang throttle by a single global polynomial, which reduces the accuracy of the recovered costates, and show the advantages of using a multiphase formulation to represent known control structures, as done for the first-stage branch, whose throttle structure is imposed through its phase sequence. Overall, the Hamiltonian remains two to three orders of magnitude smaller than its dominant terms and, at the Birkhoff-Lagrange separation interface, satisfies the branch-point transversality condition $\mathcal{H}^{(\mathrm{A})}(t_{f,\mathrm{A}})=\mathcal{H}^{(\mathrm{B})}(t_{0,\mathrm{B}})+\mathcal{H}^{(\mathrm{S})}(t_{0,\mathrm{S}})$ implied by the time linkages of Eq.~\eqref{eq:rlv-separation-linkage}, indicating that the costates recovered from the two formulations are mutually consistent.

%% file: sections/conclusions.tex
\section{Conclusions}
\label{sec:conclusions}

\par This manuscript presented theoretical and numerical advancements toward the practical use of Birkhoff \ac{ps} methods for aerospace \acp{ocp}. The Birkhoff-$a$ integration operator was constructed explicitly for \ac{flgr} grids and related to both the corresponding Lagrange differentiation matrix and, through a signed reversal, the \ac{lgr} Birkhoff-$b$ operator of the universal Birkhoff theory, from which it inherits its conditioning properties. The former relation shows that the two \ac{flgr} transcriptions share the same discrete dynamics, and that the ill-conditioning of the Lagrange operator is caused by its growing largest singular values, while the problematic smallest singular values of the Birkhoff operator are pushed away from zero by the identity block of the lifted representation. Differential constraints on nonlinear trajectory quantities were enforced directly through the \ac{ps} operators and, for control-rate limits, shown to preserve the primal and dual structure of state augmentation while avoiding its additional rate-control variables and collocation equations.

\par On the computational side, non-affine time mappings were exploited to concentrate the \ac{ps} nodes around localized trajectory features, and the composition of affine scaling and equilibration transformations provided a systematic approach to construct \ac{ocp} scaling strategies, recovering \ac{pjrn} \cite{sagliano2017_RadauPseudospectralMethod} as a particular case. All developments are implemented in \ac{raptor}, which supports heterogeneous multiphase transcription with phase-wise Lagrange or Birkhoff interpolants. The orbit-raising problem displays the expected conditioning trends of the adopted \ac{ps} operators, with the Birkhoff formulation exhibiting shorter solution times on dense grids. In the landing problem, the direct thrust-rate constraint reproduces the primal and dual solutions of state augmentation with a 31\% lower mean solution time, while the non-affine time mapping improves the coarse-grid representation of the thrust switch. In the entry problem, applying static Ruiz equilibration after \ac{pjrn} halves the mean solution time while recovering the same final latitude. Finally, the five-phase \ac{rlv} problem demonstrates the integration of these developments within a branching multiphase mission.